\documentclass[12pt]{amsart}
\usepackage{latexsym}
\usepackage{amsmath}
\usepackage{amssymb}
\usepackage{amsthm}
\usepackage[all]{xy}

\allowdisplaybreaks                          
\newtheorem{Thm}{Theorem}[section]
\newtheorem{Prop}{Proposition}[section]
\newtheorem{Cor}{Corollary}[section]
\newtheorem{Lem}{Lemma}[section]

\newtheorem{Rem}{Remark}[section]

\begin{document}

\title{Nefness of adjoint bundles for ample vector bundles of corank $3$}

\author{Andrea Luigi Tironi
}

\date{\today}



\maketitle

\begin{abstract}
Let $\mathcal{E}$ be an ample vector bundle of rank $r\geq 2$ on a
smooth complex projective variety $X$ of dimension $n$. The aim of
this paper is to describe the structure of pairs $(X,\mathcal{E})$
as above whose adjoint bundles $K_X+\det\mathcal{E}$ are not nef
for $r=n-3$. Furthermore, we give some immediate consequences of
this result in adjunction theory.
\end{abstract}

\section{Introduction}\label{sec:1}

Let $\mathcal{E}$ be an ample vector bundle of rank $r\geq 2$ on a
smooth complex projective variety $X$ of dimension $n$. We say
that a Cartier divisor $D$ on $X$ is \textit{numerically effective
(nef)} if it has non-negative intersection number with any curve
$C$ on $X$.

In the late 1980s, Ye and Zhang started to study the nefness of
the adjoint bundle $K_X+\det\mathcal{E}$, obtaining in \cite{YZ}
complete results for $r\geq n-1$.

Subsequently, Zhang investigated in \cite{Zhang} also the case
$r=n-2$, showing that in this situation $K_X+\det\mathcal{E}$ is
nef except in a few cases.

Some years later, by the adjunction theory and the remark that a
key result of \cite{W} is still valid without the extra assumption
of spannedness for $\mathcal{E}$, Maeda gave in \cite{Maeda} an
exhaustive description of all the ample vector bundles
$\mathcal{E}$ for the case $r=n-2$.

Let us note that a weaker version of the main result of
\cite{Zhang}, in which $\mathcal{E}$ is assumed to be also
generated by global sections, was already proved by
Wi$\acute{\textrm{s}}$niewski in \cite{W} and that all the above
results are a generalization of those due to Fujita \cite{Fujita}
and Ionescu \cite{I}.

Thus the purpose of this paper is to describe the structure of
pairs $(X,\mathcal{E})$ as above whose adjoint bundles
$K_X+\det\mathcal{E}$ are not nef in the next case $r=n-3$.

\smallskip

The final statement which we obtain in Section 3 is as follows:

\begin{Thm}\label{thm}
Let $\mathcal{E}$ be an ample vector bundle of rank $r\geq 2$ on a
smooth complex projective variety $X$ of dimension $n$. If
$r=n-3$, then $K_X+\det\mathcal{E}$ is nef except when
$(X,\mathcal{E})$ is one of the following:
\begin{enumerate}
\item[$(1)$] $(\mathbb{P}^n, \oplus_{i=1}^{n-3}
\mathcal{O}_{\mathbb{P}^n}(a_i)$, where all the $a_i$'s are
positive integers such that $\sum_{i=1}^{n-3}a_i\leq n$;
\item[$(2)$] $(\mathbb{P}^n, \mathcal{O}_{\mathbb{P}^n}(1)^{\oplus
n-3-a}\oplus\mathcal{V})$, where $\mathcal{V}$ is an
indecomposable Fano vector bundle on $\mathbb{P}^n$ of rank $a$
such that $3\leq a\leq n-3$, $c_1(\mathcal{V})=a+3$ and its
generic splitting types are $(3,2,1,1,...,1)$ or
$(2,2,2,1,...,1)$; moreover,
$\mathcal{V}(-1):=\mathcal{V}\otimes\mathcal{O}
(-1)$ is nef for $3\leq a\leq n-6$, but it is not globally
generated for $n\geq 6$; \item[$(3)$] $(\mathbb{Q}^n,
\oplus_{j=1}^{n-3}\mathcal{O}_{\mathbb{Q}^n}(b_j))$, where all the
$b_j$'s are positive integers such that $\sum_{j=1}^{n-3}b_j\leq
n-1$; \item[$(4)$] $(\mathbb{Q}^n,
\mathcal{O}_{\mathbb{Q}^n}(1)^{\oplus n-3-b}\oplus\mathcal{V})$,
where $\mathcal{V}$ is an indecomposable Fano vector bundle on
$\mathbb{Q}^n$ of rank $b$ such that $3\leq b\leq n-3$,
$c_1(\mathcal{V})=b+2$ and its generic splitting type is
$(2,2,1,...,1)$; moreover, $\mathcal{V}(-1)$ is nef on
$\mathbb{Q}^n$ for $3\leq b\leq n-5$, but it is not globally
generated for $n\geq 7$; \item[$(5)$] $X$ is a Fano $n$-fold of
index $n-1$ with $\mathrm{Pic}(X)$ generated by an ample line
bundle $\mathcal{H}$ and either $(\alpha)$
$\mathcal{E}\cong\mathcal{H}^{\oplus n-3}$, or $(\beta)$
$\mathcal{E}|_l\cong\mathcal{H}_l^{\oplus
n-4}\oplus\mathcal{H}_l^{\otimes 2}$ for every line $l$ of
$(X,\mathcal{H})$; \item[$(6)$] $X$ is a Fano $n$-fold of index
$n-2$ with $\mathrm{Pic}(X)$ generated by an ample line bundle
$\mathcal{L}$ and $\mathcal{E}\cong\mathcal{L}^{\oplus n-3}$;
\item[$(7)$]
$(\mathbb{P}^3\times\mathbb{P}^3,\mathcal{O}_{\mathbb{P}^3\times\mathbb{P}^3}(1,1)^{\oplus
3})$; \item[$(8)$]
$(\mathbb{P}^2\times\mathbb{Q}^3,\mathcal{O}_{\mathbb{P}^2\times\mathbb{Q}^3}(1,1)^{\oplus
2})$; \item[$(9)$] $(\mathbb{P}_{\mathbb{P}^3}(T),[\xi_T]^{\oplus
2})$, where $T$ is the tangent bundle of $\mathbb{P}^3$ and
$\xi_{T}$ is the tautological line bundle; \item[$(10)$]
$(\mathbb{P}_{\mathbb{P}^3}(\mathcal{O}(2)\oplus\mathcal{O}(1)^{\oplus
2}), [\xi]^{\oplus 2})$, $\xi$ being the tautological line bundle;
\item[$(11)$] there is a vector bundle $\mathcal{V}$ on a smooth
curve $C$ such that $X\cong\mathbb{P}_C(\mathcal{V})$; moreover,
for any fibre $F\cong\mathbb{P}^{n-1}$ of $X\to C$, we have
$\mathcal{E}|_F$ is isomorphic to one of the following:
\begin{enumerate}
\item[$(a)$] $\mathcal{O}_{\mathbb{P}^{n-1}}(1)^{\oplus n-3}$; \
$(b)$
$\mathcal{O}_{\mathbb{P}^{n-1}}(2)\oplus\mathcal{O}_{\mathbb{P}^{n-1}}(1)^{\oplus
n-4}$; \item[$(c)$] $\mathcal{O}_{\mathbb{P}^{n-1}}(2)^{\oplus
2}\oplus\mathcal{O}_{\mathbb{P}^{n-1}}(1)^{\oplus n-5}$; \ $(d)$
$\mathcal{O}_{\mathbb{P}^{n-1}}(3)\oplus\mathcal{O}_{\mathbb{P}^{n-1}}(1)^{\oplus
n-4}$;
\end{enumerate}
\item[$(12)$] $X$ is a section of a divisor of relative degree two
in a projective space $\mathbb{P}_{C}(\mathcal{G})$, where
$\mathcal{G}$ is a vector bundle of rank $n+1$ on a smooth curve
$C$; moreover, for any smooth fibre $F\cong\mathbb{Q}^{n-1}$ of
$X\to C$, where $\mathbb{Q}^{n-1}$ is a smooth quadric
hypersurface of $\mathbb{P}^n$, we have $(F,\mathcal{E}|_F)$ is
one of the following pairs:
\begin{enumerate}
\item[(i)]
$(\mathbb{Q}^{n-1},\mathcal{O}_{\mathbb{Q}^{n-1}}(1)^{\oplus
n-3})$; \item[(ii)]
$(\mathbb{Q}^{n-1},\mathcal{O}_{\mathbb{Q}^{n-1}}(2)\oplus\mathcal{O}_{\mathbb{Q}^{n-1}}(1)^{\oplus
n-4})$; \item[(iii)] $(\mathbb{Q}^4,
\mathcal{S}\otimes\mathcal{O}_{\mathbb{Q}^4}(2))$, where
$\mathcal{S}$ is the spinor bundle on
$\mathbb{Q}^4\subset\mathbb{P}^5$;
\end{enumerate}
\item[$(13)$] the map $\Phi:X\to C$ associated to the linear
system $|(n-3)K_X+(n-2)\det\mathcal{E}|$ makes $X$ a Del Pezzo
fibration over a smooth curve $ C$; moreover, any general smooth
fibre $F$ of $\Phi$ is either a Del Pezzo $(n-1)$-fold with
$\mathrm{Pic }(F)\cong\mathbb{Z}[\mathcal{O}_F(1)]$ and such that
$\mathcal{E}|_F\cong\mathcal{O}_F(1)^{\oplus n-3}$, or
$\mathbb{P}^2\times\mathbb{P}^2$ with
$\mathcal{E}|_F\cong\mathcal{O}_{\mathbb{P}^2\times\mathbb{P}^2}(1,1)^{\oplus
2}$; \item[$(14)$] there is a vector bundle $\mathcal{V}$ on a
smooth surface $S$ such that $X\cong\mathbb{P}_S(\mathcal{V})$;
moreover, for any fibre $F\cong\mathbb{P}^{n-2}$ of $X\to S$, we
have
$\mathcal{E}|_F\cong\mathcal{O}_{\mathbb{P}^{n-2}}(2)\oplus\mathcal{O}_{\mathbb{P}^{n-2}}(1)^{\oplus
n-4}$; \item[$(15)$] there exists a $\mathbb{P}^{n-2}$-fibration
$\pi:X\to S$, locally trivial in the complex (or \'{e}tale)
topology, over a smooth surface $S$ such that
$\mathcal{E}|_F\cong\mathcal{O}_{\mathbb{P}^{n-2}}(1)^{\oplus
n-3}$ for any closed fiber $F\cong\mathbb{P}^{n-2}$ of the map
$\pi$; \item[$(16)$] the map $\psi:X\to S$ associated to the
linear system $|(n-3)K_X+(n-2)\det\mathcal{E}|$ makes $X$ a
quadric fibration over a smooth surface $S$; moreover, for any
general fibre $F\cong\mathbb{Q}^{n-2}$ we have
$\mathcal{E}|_F\cong\mathcal{O}_{\mathbb{Q}^{n-2}}(1)^{\oplus
n-3}$, where $\mathbb{Q}^{n-2}$ is a smooth quadric hypersurface
of $\mathbb{P}^{n-1}$; \item[$(17)$] there is a vector bundle
$\mathcal{F}$ on a smooth $3$-fold $V$ such that
$X\cong\mathbb{P}_V(\mathcal{F})$; moreover, for any fibre
$F\cong\mathbb{P}^{n-3}$ of $X\to V$, we have
$\mathcal{E}|_F\cong\mathcal{O}_{\mathbb{P}^{n-3}}(1)^{\oplus
n-3}$; \item[$(18)$] there exists a smooth variety $X'$ and a
morphism $\varphi:X\to X'$ expressing $X$ as the blowing-up of
$X'$ at a finite set of points $B$ and an ample vector bundle
$\mathcal{E}'$ on $X'$ such that
$\mathcal{E}\otimes([\varphi^{-1}(B)])\cong\varphi^*\mathcal{E}'$
and $K_{X'}+\tau'\det\mathcal{E}'$ is nef, where $\tau' <
\frac{n-1}{n-3}$ is the nefvalue of the pair
$(X',\det\mathcal{E}')$. Moreover,
$\mathcal{E}|_E\cong\mathcal{O}_{\mathbb{P}^{n-1}}(1)^{\oplus
n-3}$ for any irreducible component $E$ of the exceptional locus
of $\varphi$; \item[$(19)$] the map $\psi :X\to X'$ associated to
the linear system $|(n-2)K_X+(n-1)\det\mathcal{E}|$ is a
birational morphism which contracts an extremal face spanned by
extremal rays $R_i$ for some $i$ in a finite set of index. Let
$\psi_i:X\to X_i$ be the contraction associated to $R_i$. Then
each $\psi_i$ is birational and of divisorial type; moreover, if
$E_i$ is an exceptional divisor of $\psi_i$, then
$(E_i,[E_i]_{E_i},\mathcal{E}|_{E_i})\cong
(\mathbb{P}^{n-1},\mathcal{O}_{\mathbb{P}^{n-1}}(-2),\mathcal{O}_{\mathbb{P}^{n-1}}(1)^{\oplus
n-3})$; \item[$(20)$] the map $\phi:X\to X'$ associated to the
linear system $|(n-3)K_X+(n-2)\det\mathcal{E}|$ is a birational
morphism which contracts an extremal face. Let $R_i$ be the
extremal rays spanning this face for some $i$ in a finite set of
index. Call $\rho_i:X\to W_i$ the contraction associated to one of
the $R_i$. Then we have each $\rho_i$ is birational and of
divisorial type; if $D_i$ is one of the exceptional divisors and
$Z_i=\rho_i(D_i)$, we have $\mathrm{dim }Z_i\leq 1$ and one of the
following possibilities can occur:
\begin{enumerate}
\item[(j)] $\mathrm{dim } Z_i=0, D_i\cong\mathbb{P}^{n-1}$ and
$([D_i]_{D_i},\mathcal{E}|_{D_i})\cong
(\mathcal{O}_{\mathbb{P}^{n-1}}(-1),\mathcal{O}_{\mathbb{P}^{n-1}}(1)^{\oplus
n-4}\oplus\mathcal{O}_{\mathbb{P}^{n-1}}(2))$; \item[(jj)]
$\mathrm{dim } Z_i=0, D_i$ is a (possible singular) quadric
$\mathbb{Q}^{n-1}\subset\mathbb{P}^{n}$ and
$[D_i]_{D_i}=\mathcal{O}_{\mathbb{Q}^{n-1}}(-1)$; moreover,
$\mathcal{E}|_{D_i}\cong\mathcal{O}_{\mathbb{Q}^{n-1}}(1)^{\oplus
n-3}$; \item[(jjj)] $\mathrm{dim } Z_i=1$, $W_i$ and $Z_i$ are
smooth projective varieties and $\rho_i$ is the blow-up of $W_i$
along $Z_i$; moreover,
$\mathcal{E}|_{F_i}\cong\mathcal{O}_{F_i}(1)^{\oplus n-3}$ and
${\mathcal{O}_{D_i}(D_i)}|_{F_i}=\mathcal{O}_{\mathbb{P}^{n-2}}(-1)$
for any fibre $F_i\cong\mathbb{P}^{n-2}$ of $D_i\to Z_i$.
\end{enumerate}
Moreover, the map $\phi$ is a composition of disjoint extremal
contractions as in $\mathrm{(j), (jj)}$ and $\mathrm{(jjj)}$.
\end{enumerate}
\end{Thm}

\smallskip

\noindent The Theorem \ref{thm} is a generalization of a result
due to Fujita \cite[theorem 4]{Fujita} and its proof makes use
principally of \cite{AM} and all the above papers \cite{Maeda},
\cite{YZ} and \cite{Zhang}, together with some recent results
(\cite{AN} and \cite{Ohno}) about the nef value of pairs
$(X,\det\mathcal{E})$.

Note that Theorem \ref{thm} is complete for $n=5$ (or
$\mathrm{Pic}(X)\neq\mathbb{Z}$). When $n\geq 6$ and $X$ is a Fano
manifold of index $i(X)\geq n-2$ with
$\mathrm{Pic}(X)=\mathbb{Z}[\mathcal{O}_X(1)]$, a finest
description of $\mathcal{E}$ is equivalent to have in this setting
a complete classification of ample vector bundles
$\mathcal{F}:=\mathcal{E}\oplus\mathcal{O}_X(1)$ of rank $n-2\geq
4$ such that either $c_1(\mathcal{F})=-K_X$ (see also
\cite[theorem 5.1 (2)(i)]{AM}), or $\mathbb{P}_X(\mathcal{F}):=W$
is a Fano manifold of index $i(W)=\frac{\mathrm{dim } W-1}{2}\geq
4$ (the case $i(W)=3$ implies $n=5$; see \cite[theorem 1.3]{NO}
for the case $i(W)=2$). However, a more satisfactory and effective
description of $\mathcal{E}$ when
$\mathrm{Pic}(X)=\mathbb{Z}[\mathcal{O}_X(1)]$ is given in
Corollary \ref{corollary thm} under the extra assumption that
$\mathcal{E}\otimes\mathcal{O}_X(-1)$ is globally generated on
$X$.

In the last section, we consider two immediate applications of
Theorem \ref{thm}. More precisely, in $\S$\ref{3.1} we extend
\cite[theorem 3]{LM},
by giving a description of ample vector bundles $\mathcal{E}$ of
rank $r\geq 2$ on a smooth projective variety $X$ of dimension $n$
which admit a global section $s\in\Gamma (\mathcal{E})$ whose zero
locus $Z:=(s)_0\subset X$ is a smooth subvariety of the expected
dimension $n-r\geq 3$ and such that $K_Z+(\mathrm{dim }Z-3)H_Z$ is
not nef (Propositions \ref{dim Z=3} and \ref{dim Z>3}), where
$H_Z$ denotes the restriction to $Z$ of an ample line bundle $H$
on $X$. These results partially overlap with those obtained in
\cite[$\S\S 4,5,6$ and $7$]{AO}, but they are proved in a simple
and different way, and together with \cite[theorems $1,2$ and
$3$]{LM}, they describe the structure of all pairs
$(X,\mathcal{E})$ as above such that $(Z,H_Z)$ is a special
polarized variety until the second reduction map in the sense of
the adjunction theory (e.g., see \cite[chapter 7]{BS} for a
complete classification of such pairs $(Z,H_Z)$). Finally, in
$\S$\ref{3.2} we study classical scrolls over manifolds of
dimension five and Picard number one (Corollary \ref{corollary}
and Remark \ref{remark}), extending some results of \cite{BS2} and
\cite[$\S 3$]{T} about classical scrolls which are not
adjunction-theoretic scrolls (we refer to $\S 1$ and the first
part of $\S$\ref{3.2} for the difference between the two
definitions of scrolls).

\section{Notation and background material.}

In this note varieties are always assumed to be defined over the
complex number field $\mathbb{C}$. We use the standard notation
from algebraic geometry. The words ``vector bundles'' and
``locally free sheaves'' are used interchangeably. Let $X$ be a
smooth irreducible projective variety of dimension $n$ (for
simplicity, $n$-fold). The group of line bundles on $X$ is denoted
by Pic$(X)$. Moreover, we denote by $\rho(X)$ the Picard number of
$X$. The tensor products of line bundles are denoted additively.
The pull-back $j^*\mathcal{E}$ of a vector bundle $\mathcal{E}$ on
$X$ by an embedding $j:V\to X$ is denoted by $\mathcal{E}|_V$ (or,
for simplicity, by $\mathcal{E}_V$ when no confusion arises). The
canonical bundle of an $n$-fold $X$ is denoted by $K_X$. An
$n$-fold $X$ is said to be a \textit{Fano manifold} if its
anticanonical bundle $-K_X$ is ample. For Fano manifolds $X$, the
largest integer $i(X)$ which divides $-K_X$ in Pic$(X)$ is called
the \textit{index} of $X$.

A \textit{polarized $n$-fold} is a pair $(X,L)$ consisting of an
$n$-fold $X$ and an ample line bundle $L$ on $X$. A polarized
$n$-fold $(X,L)$ is said to be a \textit{classical scroll} over a
smooth variety $Y$ if $(X,L)\cong
(\mathbb{P}_Y(\mathcal{V}),\xi_{\mathcal{V}})$ for some ample
vector bundle $\mathcal{V}$ on $Y$, where $\xi_{\mathcal{V}}$ is
the tautological line bundle on the projective space
$\mathbb{P}_Y(\mathcal{V})$ associated to $\mathcal{V}$. We say
that $(X,L)$ is an \textit{adjunction--theoretic scroll}
(respectively a \textit{quadric fibration}, respectively a
\textit{Del Pezzo fibration}, respectively a \textit{Mukai
fibration}) over a normal variety $Y$ of dimension $m$ if there
exists a surjective morphism with connected fibres $p:X\to Y$ and
an ample line bundle $H$ on $Y$, such that $K_X+(n-m+1)L\simeq
p^*H$ (respectively $K_X+(n-m)L\simeq p^*H$, respectively
$K_X+(n-m-1)L\simeq p^*H$, respectively $K_X+(n-m-2)L\simeq
p^*H$). We say that $(X,L)$ is a \textit{Del Pezzo variety}
(respectively a \textit{Mukai variety}) if $K_X\simeq -(n-1)L$
(respectively $K_X\simeq -(n-2)L$). Moreover, for general results
about the adjunction theory, we refer to \cite{BS} and
\cite{Sommese}.

A part of Mori's theory of extremal rays is to be used throughout
the paper. So, let $Z_1(X)$ be the free abelian group generated by
integral curves on an $n$-fold $X$. The intersection pairing gives
a bilinear map Pic$(X)\times Z_1(X)\to\mathbb{Z}$ and the
\textit{numerical equivalence} $\equiv$ is defined so that the
pairing
$$((\mathrm{Pic}(X)/\equiv)\otimes\mathbb{Q})\times
((Z_1(X)/\equiv)\otimes\mathbb{Q})\to\mathbb{Q}$$ is
non-degenerate. The \textit{closed cone of curves}
$\overline{NE}(X)$ is the closed convex cone generated by
effective $1$-cycles in the $\mathbb{R}$-vector space
$(Z_1(X)/\equiv)\otimes\mathbb{R}$. We say that
$L\in\mathrm{Pic}(X)$ is \textit{nef} if the numerical class of
$L$ in $(\mathrm{Pic}(X)/\equiv)\otimes\mathbb{R}$ gives a
non-negative function on $\overline{NE}(X)-\{0\}$. Let $Z$ be a
$1$-cycle on $X$. We denote by $[Z]$ the numerical class of $Z$ in
$(Z_1(X)/\equiv)\otimes\mathbb{R}$. A half line
$R=\mathbb{R}_+[Z]$ in $\overline{NE}(X)$ is called an
\textit{extremal ray} if

\smallskip

(1) $K_X\cdot Z<0$, and

\smallskip

(2) if $z_1,z_2\in\overline{NE}(X)$ satisfy $z_1+z_2\in R$, then
$z_1,z_2\in R$.

\smallskip

\noindent A rational (possibly singular) reduced and irreducible
curve $C$ on $X$ is called an \textit{extremal rational curve} if
$\mathbb{R}_+[C]$ is an extremal ray and $(-K_X)\cdot C\leq n+1$.
Let
$$\overline{NE}(X)^{+}=\{z\in\overline{NE}(X) \ | \ K_X\cdot z\geq
0\}.$$

\noindent Then we have the following basic theorem in the Mori
theory.

\begin{Thm}[Cone Theorem]\label{cone theorem}
Let $X$ be a smooth projective variety. Then $\overline{NE}(X)$ is
the smallest closed convex cone containing $\overline{NE}(X)^+$
and all the extremal rays:
$$\overline{NE}(X)=\overline{NE}(X)^{+}+\sum_i R_i,$$ where the
$R_i$ are extremal rays of $\overline{NE}(X)$ for $X$. For any
open convex cone $V$ containing $\overline{NE}(X)^{+}-\{0\}$ there
exist only a finite number of extremal rays that do not lie in
$V\cup\{0\}$. Furthermore, every extremal ray is spanned by a
numerical class of an extremal rational curve.
\end{Thm}

For the proof of the above result, we refer to \cite[theorem
1.5]{Mori} and \cite[(1.2)]{Mori1}. In the next section, we will
use also the following well-known

\begin{Lem}\label{lemma}
Let $\mathcal{E}$ be an ample vector bundle of rank $r$ on a
rational curve $C$. Then $\det\mathcal{E}\cdot C\geq r$.
\end{Lem}

Now, let us give here some technical results about vector bundles
on some Fano $n$-folds $X$ with $\rho(X)=1$.

\begin{Lem}\label{lemma 1}
Let $\mathcal{V}$ be an ample vector bundle of rank $r$ on a Fano
$n$-fold $X$ with $\rho(X)=1$. Assume that $c_1(\mathcal{V})\leq
i(X)-1$, where $i(X)$ is the index of $X$. If
$\frac{n-r+1}{i(X)-c_1(\mathcal{V})}<2$ then
$\mathcal{V}\otimes\mathcal{O}_X(-1)$ is a nef vector bundle on
$X$.
\end{Lem}

\noindent\textbf{Proof.} By adjunction we have
$-K_{\mathbb{P}(\mathcal{V})}\simeq
r\xi_{\mathcal{V}}+\pi^*\mathcal{O}_X(i(X)-c_1(\mathcal{V}))$,
where $\xi_{\mathcal{V}}$ is the tautological line bundle on
$\mathbb{P}(\mathcal{V})$ and $\pi:\mathbb{P}(\mathcal{V})\to X$
is the projection map. Since $i(X)-c_1(\mathcal{V})>0$ and
$\xi_{\mathcal{V}}$ is ample, we see that
$\mathbb{P}(\mathcal{V})$ is a Fano $(n+r-1)$-fold. Let $R$ be an
extremal ray of $\mathbb{P}(\mathcal{V})$ different from a line in
a fiber of $\pi$. Then we get
$$n+1\geq -K_{\mathbb{P}(\mathcal{V})}\cdot R
\geq r+[i(X)-c_1(\mathcal{V})]\pi^*\mathcal{O}_X(1)\cdot R.$$ This
gives
$2>\frac{n-r+1}{i(X)-c_1(\mathcal{V})}\geq\pi^*\mathcal{O}_X(1)\cdot
R,$ showing that $0\leq \pi^*\mathcal{O}_X(1)\cdot R'\leq 1$ for
any extremal ray $R'$ on $\mathbb{P}(\mathcal{V})$. So, let $C$ be
an effective irreducible curve on $\mathbb{P}(\mathcal{V})$. From
Theorem \ref{cone theorem} it follows that $C\equiv
\sum_{i}\alpha_i R_i$ with $\alpha_i>0$, where each $R_i$ is an
extremal ray of $\mathbb{P}(\mathcal{V})$. Therefore we deduce
that
$$[\xi_{\mathcal{V}}-\pi^*\mathcal{O}_X(1)]\cdot
C=\sum_i\alpha_i[\xi_{\mathcal{V}}-\pi^*\mathcal{O}_X(1)]\cdot
R_i\geq 0.$$ So by definition
$\xi_{\mathcal{V}}-\pi^*\mathcal{O}_X(1)$ is nef on
$\mathbb{P}(\mathcal{V})$, i.e.
$\mathcal{V}\otimes\mathcal{O}_X(-1)$ is nef on $X$. $\square$

\begin{Lem}\label{lemma 2}
Let $\mathcal{V}$ be a globally generated vector bundle of rank
$r<n$ on $\mathbb{P}^n$ with
$c_1:=c_1(\mathcal{V})\leq\frac{n}{2}$. If $r\geq c_1+1$, then one
of the following possibilities can occurs:
\begin{enumerate}
\item[(1)]
$\mathcal{V}\cong\mathcal{O}_{\mathbb{P}^n}(c_1)\oplus\mathcal{O}_{\mathbb{P}^n}^{\oplus
r-1}$; \item[(2)]
$\mathcal{V}\cong\mathcal{O}_{\mathbb{P}^n}(1)^{\oplus
c_1}\oplus\mathcal{O}_{\mathbb{P}^n}^{\oplus r-c_1}$; \item[(3)]
$\mathcal{V}\cong\mathcal{O}_{\mathbb{P}^n}(2)\oplus\mathcal{O}_{\mathbb{P}^n}(1)^{\oplus
c_1-2}\oplus\mathcal{O}_{\mathbb{P}^n}^{\oplus r-c_1+1}$;
\item[(4)]
$\mathcal{V}\cong\mathcal{O}_{\mathbb{P}^n}(3)\oplus\mathcal{O}_{\mathbb{P}^n}(1)^{\oplus
c_1-3}\oplus\mathcal{O}_{\mathbb{P}^n}^{\oplus r-c_1+2}$;
\item[(5)] $\mathcal{V}\cong\mathcal{O}_{\mathbb{P}^n}(2)^{\oplus
2}\oplus\mathcal{O}_{\mathbb{P}^n}(1)^{\oplus
c_1-4}\oplus\mathcal{O}_{\mathbb{P}^n}^{\oplus r-c_1+2}$;
\item[(6)] $c_1\geq 5$ and there exist the following exact
sequences:
$$0\to\mathcal{O}_{\mathbb{P}^n}\to\mathcal{V}\to\mathcal{V}^{r-1}\to
0, \ ... \ ,
0\to\mathcal{O}_{\mathbb{P}^n}\to\mathcal{V}^{k+1}\to\mathcal{V}^{k}\to
0,$$ where $2\leq k\leq c_1-3$, $\mathrm{rk}\mathcal{V}^i=i$,
$c_1(\mathcal{V}^i)=c_1$, all the $\mathcal{V}^i$'s are globally
generated and $(s_k)_0\neq\emptyset$ for a generic section $s_k$
of $\mathcal{V}^k$.
\end{enumerate}
\end{Lem}

\noindent\textbf{Proof.} Let $Z_r=(s_r)_0$ be the zero locus of a
generic section $s_r$ of $\mathcal{V}$. If $Z_r\neq\emptyset$,
then $Z_r$ is smooth and of dimension $n-r$. By adjunction
$K_{Z_r}=\mathcal{O}_{\mathbb{P}^n}(-n-1+c_1)|_{Z_r}$. Thus
$n+1-c_1\leq\dim Z_r+1=n-r+1$, i.e. $r\leq c_1$, but this is a
contradiction. Hence $Z_r=\emptyset$ and by induction (see, e.g.,
\cite[(4.3.2), p.83]{OSS}) we obtain the following exact
sequences:
$$0\to\mathcal{O}_{\mathbb{P}^n}\to\mathcal{V}\to\mathcal{V}^{r-1}\to
0, \ ... \ ,
0\to\mathcal{O}_{\mathbb{P}^n}\to\mathcal{V}^{k+1}\to\mathcal{V}^{k}\to
0,$$ where $k\leq c_1$, $\mathrm{rk}\mathcal{V}^i=i$,
$c_1(\mathcal{V}^i)=c_1$, all the $\mathcal{V}^i$'s are globally
generated and $Z_k=(s_k)_0\neq\emptyset$ for a generic section
$s_k$ of $\mathcal{V}^k$. If $k=1$ then
$\mathcal{V}^1=\mathcal{O}_{\mathbb{P}^n}(c_1)$ and by induction
we deduce that $\mathcal{V}\cong\mathcal{O}_{\mathbb{P}^n}^{\oplus
r-1}\oplus\mathcal{O}_{\mathbb{P}^n}(c_1)$. So $2\leq k\leq c_1$.
Assume that $k=c_1\geq 2$. Since $Z_{c_1}=Z_k\neq\emptyset$,
$Z_{c_1}$ is smooth and of dimension $n-c_1$. Moreover, since
$2(n-c_1)-n=n-2c_1\geq 0$, we obtain that $Z_{c_1}$ is also
irreducible. Thus by the Kobayashi-Ochiai Theorem we have
$Z_{c_1}=\mathbb{P}^{n-c_1}$. Note that
$\mathcal{N}_{Z_{c_1}/\mathbb{P}^n}\cong\mathcal{O}_{Z_{c_1}}(1)^{\oplus
c_1}$. Since
$$\mathcal{V}^{c_1}_{|Z_{c_1}}\cong\mathcal{N}_{Z_{c_1}/\mathbb{P}^n}\cong\mathcal{O}_{\mathbb{P}^{n-c_1}}(1)^{\oplus
c_1}$$ and $\dim Z_{c_1}=n-c_1\geq c_1\geq 2$, by \cite[Ch.I
(2.3.2)]{OSS} we conclude that $\mathcal{V}^{c_1}$ splits, and
then
$\mathcal{V}\cong\mathcal{O}_{\mathbb{P}^n}^{r-c_1}\oplus\mathcal{O}_{\mathbb{P}^n}(1)^{\oplus
c_1}$. Suppose now that $k=c_1-1\geq 2$. By arguing as above, we
have $Z_{c_1-1}\neq\emptyset$ is smooth, irreducible and of
dimension $n-c_1+1$. Therefore by the Kobayashi-Ochiai Theorem we
know that
$Z_{c_1-1}\cong\mathbb{Q}^{n-c_1+1}\subset\mathbb{P}^{n}$. Hence
$\mathcal{N}_{Z_{c_1-1}/\mathbb{P}^n}\cong\mathcal{O}_{\mathbb{P}^n}(2)_{|Z_{c_1-1}}\oplus\mathcal{O}_{\mathbb{P}^n}(1)^{c_1-2}_{|Z_{c_1-1}}\cong
\mathcal{V}^{c_1-1}_{|Z_{c_1-1}}$. Since $\dim
Z_{c_1-1}=n-c_1+1\geq c_1+1=k+2\geq 4$, we see that there exists a
linear $\mathbb{P}^2\subset Z_{c_1-1}$ such that
$\mathcal{V}^{c_1-1}_{|\mathbb{P}^2}$ splits. This shows that
$\mathcal{V}\cong\mathcal{O}_{\mathbb{P}^n}^{r-c_1+1}\oplus\mathcal{O}_{\mathbb{P}^n}(2)\oplus\mathcal{O}_{\mathbb{P}^n}(1)^{\oplus
c_1-2}$. Therefore $2\leq k\leq c_1-2$. Finally, suppose that
$k=c_1-2\geq 2$. By arguing as in the above cases, we obtain that
$Z_{c_1-2}\neq\emptyset$ is a smooth irreducible Del Pezzo
$(n-c_1+2))$-fold with index
$i(Z_{c_1-2})=n-c_1+1\geq\frac{n}{2}+1>\frac{\dim
Z_{c_1-2}}{2}+1$. Then by \cite{W1} we know that
$\rho(Z_{c_1-2})=1$. Furthermore, since $\dim
Z_{c_1-2}=n-c_1+2\geq\frac{n}{2}+2\geq c_1+2\geq 6$, by the
classification of Del Pezzo manifolds (see, e.g., \cite{F2}), we
conclude that $Z_{c_1-2}$ is one of the following manifolds:
\begin{enumerate}
\item[(a)] a cubic hypersurface in $\mathbb{P}^{n-c_1+3}$;
\item[(b)] a complete intersection of two quadric hypersurfaces
$\mathbb{Q}_i\subset\mathbb{P}^{n-c_1+4}$; \item[(c)]
$\mathbb{G}(1,4)\subset\mathbb{P}^9$.
\end{enumerate}
Note that Case (c) does not occur, since in this situation
$n=c_1+4\leq\frac{n}{2}+4$, i.e.
$\mathbb{G}(1,4)\subset\mathbb{P}^n$ with $n\leq 8$, a
contradiction. On the other hand, in Case (a) we get
$$\mathcal{V}^{c_1-2}_{|Z_{c_1-2}}\cong\mathcal{N}_{Z_{c_1-2}/\mathbb{P}^n}\cong\mathcal{O}_{\mathbb{P}^n}(3)_{|Z_{c_1-2}}
\oplus\mathcal{O}_{\mathbb{P}^n}(1)^{\oplus c_1-3}_{|Z_{c_1-2}},$$
while in Case (b) we have
$$\mathcal{V}^{c_1-2}_{|Z_{c_1-2}}\cong\mathcal{N}_{Z_{c_1-2}/\mathbb{P}^n}\cong\mathcal{O}_{\mathbb{P}^n}(2)^{\oplus
2}_{|Z_{c_1-2}}\oplus\mathcal{O}_{\mathbb{P}^n}(1)^{\oplus
c_1-4}_{|Z_{c_1-2}}.$$ By \cite{Pr} we see that in both cases
there exists a linear $\mathbb{P}^2\subset Z_{c_1-2}$. So
$\mathcal{V}^{c_1-2}$ splits and by induction we obtain that
$\mathcal{V}^r$ is either
$\mathcal{O}_{\mathbb{P}^n}(3)\oplus\mathcal{O}_{\mathbb{P}^n}(1)^{\oplus
c_1-3}\oplus\mathcal{O}_{\mathbb{P}^n}^{\oplus r-c_1+2}$, or
$\mathcal{O}_{\mathbb{P}^n}(2)^{\oplus
2}\oplus\mathcal{O}_{\mathbb{P}^n}(1)^{\oplus
c_1-4}\oplus\mathcal{O}_{\mathbb{P}^n}^{\oplus r-c_1+2}$.
$\square$

\begin{Cor}\label{corollary 3}
Let $\mathcal{V}$ be a globally generated vector bundle of rank
$r<n$ on $\mathbb{P}^n$ with $c_1(\mathcal{V})=3$. If $n\geq 6$
then $\mathcal{V}$ splits.
\end{Cor}

\noindent\textbf{Proof.} If $r\geq 4$ then we conclude by Lemma
\ref{lemma 2}. Since $\mathcal{V}(1)$ is ample and
$c_1(\mathcal{V}(1))=c_1(\mathcal{V})+r=3+r$, if $r\leq 2$ then by
\cite[theorem (9.1)]{APW} we know that $\mathcal{V}(1)$ splits,
i.e. $\mathcal{V}$ splits. So let $r=3$. Let $s_3$ be a general
section of $\mathcal{V}$ and put $Z_3=(s_3)_0$. If $Z_3=\emptyset$
then we obtain the following exact sequence
$$0\to\mathcal{O}_{\mathbb{P}^n}\to\mathcal{V}\to\mathcal{V}^2\to
0,$$ where $\mathcal{V}^2$ is a rank-2 vector bundle on
$\mathbb{P}^n$. Since $\mathrm{rk}\mathcal{V}^2=2$ and $n\geq 6$,
from \cite[(9.1)]{APW} we deduce by a similar argument as above
that $\mathcal{V}^2$ splits, i.e. $\mathcal{V}$ splits. Therefore
we can assume that $Z_3\neq\emptyset$. Then $Z_3$ is smooth and of
dimension $n-3$. Since $2(n-3)-n=n-6\geq 0$, we see that $Z_3$ is
irreducible and by adjunction
$K_{Z_3}=(-n-1+c_1)\mathcal{O}_{\mathbb{P}^n}(1)|_{Z_3}=(-n+2)\mathcal{O}_{\mathbb{P}^n}(1)|_{Z_3}$.
Since $n-2=\dim Z_3+1$, by the Kobayashi-Ochiai Theorem we obtain
that $Z_3=\mathbb{P}^{n-3}\subset\mathbb{P}^n$. Hence
$$\mathcal{V}|_{Z_3}\cong\mathcal{N}_{Z_3/\mathbb{P}^n}\cong\mathcal{O}_{\mathbb{P}^{n-3}}(1)^{\oplus
3}$$ with $\dim Z_3=n-3\geq 3$, and this shows that $\mathcal{V}$
splits. $\square$

\begin{Cor}\label{corollary 4}
Let $\mathcal{V}$ be a globally generated vector bundle of rank
$r<n$ on $\mathbb{P}^n$ with $c_1(\mathcal{V})=4$. If $n\geq 8$
then $\mathcal{V}$ splits.
\end{Cor}

\noindent\textbf{Proof.} If $r\geq 5$ then we obtain the statement
by Lemma \ref{lemma 2}. Let $2\leq r\leq 4$. Therefore we have
either (i) $Z_r:=(s_r)_0\neq\emptyset$ for a generic section $s_r$
of $\mathcal{V}$, or by induction (ii) there exist the following
exact sequences
$$0\to\mathcal{O}_{\mathbb{P}^n}\to\mathcal{V}\to\mathcal{V}^{r-1}\to
0, \ ... \ ,
0\to\mathcal{O}_{\mathbb{P}^n}\to\mathcal{V}^{k+1}\to\mathcal{V}^{k}\to
0,$$ where $1\leq k\leq r-1\leq 3$, $\mathrm{rk}\mathcal{V}^i=i$,
$c_1(\mathcal{V}^i)=c_1$, all the $\mathcal{V}^i$'s are globally
generated and $Z_k:=(s_k)_0\neq\emptyset$ for a generic section
$s_k$ of $\mathcal{V}^k$. In Case (i), $Z_r\neq\emptyset$ is a
smooth irreducible $(n-r)$-fold with
$K_{Z_r}+(n-3)\mathcal{O}_{\mathbb{P}^n}(1)|_{Z_r}=\mathcal{O}_{\mathbb{P}^n}$
and $n-3=\dim Z_r+r-3$. Since $2\leq r\leq 4$, we get one of the
following possibilities:
\begin{enumerate}
\item[$(\alpha)$] $r=4$ and $Z_4=\mathbb{P}^{n-4}$;
\item[$(\beta)$] $r=3$ and
$Z_3=\mathbb{Q}^{n-3}\subset\mathbb{P}^{n-2}$; \item[$(\gamma)$]
$r=2$.
\end{enumerate}
In Case $(\gamma)$, since $\mathcal{V}(1)$ is ample and
$c_1(\mathcal{V}(1))=c_1(\mathcal{V})+r=6$, from \cite[(9.1)]{APW}
it follows that $\mathcal{V}$ splits. Moreover, by similar
arguments as in the proof of Corollary \ref{corollary 3}, we see
that $\mathcal{V}$ splits also in Cases $(\alpha)$ and $(\beta)$.

\smallskip

\noindent Consider now Case (ii). Note that $k=3$, since otherwise
$k\leq 2$ and we could deduce by \cite[(9.1)]{APW} that
$\mathcal{V}^k$ splits, i.e. $\mathcal{V}$ splits. Thus $Z_3\neq
\emptyset$ is a smooth irreducible $(n-3)$-fold. Moreover, by
adjunction $K_{Z_3}=\mathcal{O}_{\mathbb{P}^n}(n-3)|_{Z_3}$ and
from the Kobayashi-Ochiai Theorem we get
$Z_3=\mathbb{Q}^{n-3}\subset\mathbb{P}^n$. Thus
$\mathcal{V}|_{Z_3}\cong\mathcal{N}_{Z_3/\mathbb{P}^n}\cong\mathcal{O}_{\mathbb{Q}^{n-3}}(2)\oplus\mathcal{O}_{\mathbb{Q}^{n-3}}(1)^{\oplus
2}$ and since $\dim Z_3=n-3\geq 5$, we conclude that there exists
a $\mathbb{P}^2\subset Z_3$, i.e. $\mathcal{V}$ splits. $\square$

\bigskip

Finally, let us give here also some useful results about globally
generated vector bundles on a smooth quadric hypersurface
$\mathbb{Q}^n\subset\mathbb{P}^{n+1}$.

\begin{Lem}\label{lemma 3}
Let $\mathcal{V}$ be a globally generated vector bundle of rank
$r<n$ on $\mathbb{Q}^n\subset\mathbb{P}^{n+1}$ with
$c_1:=c_1(\mathcal{V})\leq\frac{n-3}{2}$. If $r\geq c_1+2$, then
either $\mathcal{V}$ splits, or $c_1\geq 2$ and there exist exact
sequences
$$0\to\mathcal{O}_{\mathbb{Q}^n}\to\mathcal{V}\to\mathcal{V}^{r-1}\to
0, \ ... \ ,
0\to\mathcal{O}_{\mathbb{Q}^n}\to\mathcal{V}^{k+1}\to\mathcal{V}^{k}\to
0,$$ where $2\leq k\leq c_1$, $\mathrm{rk}\mathcal{V}^i=i$,
$c_1(\mathcal{V}^i)=c_1$, all the $\mathcal{V}^i$'s are globally
generated and $(s_k)_0\neq\emptyset$ is an irreducible smooth
$(n-k)$-fold for a generic section $s_k$ of $\mathcal{V}^k$.
\end{Lem}

\noindent\textbf{Proof.} Let $Z_r$ be the zero locus of a general
section $s_r$ of $\mathcal{V}$. If $Z_r\neq\emptyset$ then $Z_r$
is smooth and of dimension $n-r$. Since by adjunction
$K_{Z_r}+(n-c_1)\mathcal{O}_{\mathbb{Q}^n}(1)|_{Z_r}=\mathcal{O}_{Z_r}$,
we deduce that $n-c_1\leq\dim Z_r+1=n-r+1$, i.e. $r\leq c_1+1$, a
contradiction. Thus $Z_r=\emptyset$ and by induction there exist
the following exact sequences:
$$0\to\mathcal{O}_{\mathbb{Q}^n}\to\mathcal{V}\to\mathcal{V}^{r-1}\to
0, \ ... \ ,
0\to\mathcal{O}_{\mathbb{Q}^n}\to\mathcal{V}^{k+1}\to\mathcal{V}^{k}\to
0,$$ where $1\leq k\leq c_1+1$, $\mathrm{rk}\mathcal{V}^i=i$,
$c_1(\mathcal{V}^i)=c_1$, all the $\mathcal{V}^i$'s are globally
generated and $(s_k)_0\neq\emptyset$ is smooth and of dimension
$(n-k)$ for a generic section $s_k$ of $\mathcal{V}^k$. Note that
$\mathcal{V}$ splits for $k=1$. Thus suppose that $2\leq k\leq
c_1+1$. Moreover, if $c_1\leq 1$, then from \cite{PSW1} we know
that $\mathcal{V}$ splits again. So we can assume also that
$c_1\geq 2$. If $k=c_1+1\geq 3$, then $Z_{c_1+1}=Z_k\neq\emptyset$
is smooth and of dimension $n-c_1-1$. Since
$$2(n-k)-n-1=n-2k-1=n-2(c_1+1)-1=n-2c_1-3\geq 0,$$ we see that
$Z_{c_1+1}$ is also irreducible. Furthermore, by adjunction we get
$$K_{Z_{c_1+1}}+(n-c_1)\mathcal{O}_{\mathbb{Q}^n}(1)|_{Z_{c_1+1}}\simeq\mathcal{O}_{Z_{c_1+1}},$$
where $n-c_1=\dim Z_{c_1+1}+1$. From the Kobayashi-Ochiai Theorem
it follows that
$Z_{c_1+1}=\mathbb{P}^{n-c_1-1}\subset\mathbb{Q}^n$. Hence
$n-c_1-1\leq\frac{n}{2}$, i.e. $n\leq 2c_1+2$, but this is absurd.
$\square$

\begin{Cor}\label{corollary 5}
Let $\mathcal{V}$ be a globally generated vector bundle of rank
$r<n$ on $\mathbb{Q}^n\subset\mathbb{P}^{n+1}$ with
$c_1(\mathcal{V})=2$. If $n\geq 7$ then $\mathcal{V}$ splits.
\end{Cor}

\noindent\textbf{Proof.} Note that
$\mathcal{V}':=\mathcal{V}\otimes\mathcal{O}_{\mathbb{Q}^n}(1)$ is
ample with
$c_1(\mathcal{V}')=c_1(\mathcal{V})+\mathrm{rk}\mathcal{V}=2+\mathrm{rk}\mathcal{V}$.
Since $n\geq 7$, if $r=2$ then $\mathcal{V}'$ is a Fano bundle of
rank-2 on $\mathbb{Q}^n$ and by \cite[(2.4)(2)]{APW} we can
conclude that $\mathcal{V}'$ splits, i.e. $\mathcal{V}$ splits.
Let $r\geq 3$. If $r\geq 4$, then from Lemma \ref{lemma 3} we
deduce that either $\mathcal{V}$ splits, or there exist exact
sequences
$$0\to\mathcal{O}_{\mathbb{Q}^n}\to\mathcal{V}\to\mathcal{V}^{r-1}\to
0, \ ... \ ,
0\to\mathcal{O}_{\mathbb{Q}^n}\to\mathcal{V}^{3}\to\mathcal{V}^{2}\to
0,$$ where $\mathrm{rk}\mathcal{V}^i=i$, $c_1(\mathcal{V}^i)=2$,
all the $\mathcal{V}^i$'s are globally generated and
$(s_2)_0\neq\emptyset$ is smooth and of dimension $(n-2)$ for a
general section $s_2$ of $\mathcal{V}^2$. Since $\mathcal{V}^2(1)$
is a Fano bundle of rank $2$ on $\mathbb{Q}^n$ with $n\geq 7$, by
\cite[main theorem (2.4)(2)]{APW} we see that $\mathcal{V}^2(1)$
splits, i.e. $\mathcal{V}^2$ splits. Therefore, from the above
exact sequences we deduce that also $\mathcal{V}$ splits.

Finally, suppose that $r=3$. Let $Z_3=(s_3)_0$ be the zero locus
of a general section $s_3$ of $\mathcal{V}$. If $Z_3=\emptyset$
then there exists an exact sequence
$0\to\mathcal{O}_{\mathbb{Q}^n}\to\mathcal{V}\to\mathcal{V}^{2}\to
0$, where $\mathcal{V}^2$ splits since $n\geq 7$, that is,
$\mathcal{V}$ splits. Hence $Z_3\neq\emptyset$ is smooth and of
dimension $n-3$. Since $2(n-3)-n-1=n-7\geq 0$, we see that $Z_3$
is also irreducible. By adjunction we obtain that
$K_{Z_3}+(n-2)\mathcal{O}_{\mathbb{Q}^n}(1)|_{Z_3}=\mathcal{O}_{Z_3}$
with $n-2=\dim Z_3+1$. Therefore
$Z_3=\mathbb{P}^{n-3}\subset\mathbb{Q}^n$ and this implies that
$n-3\leq\frac{n}{2}$, i.e. $n\leq 6$, but this gives a
contradiction. $\square$

\bigskip

Let $X$ be an $n$-fold. If $R$ is an extremal ray, then its
\textit{length} $l(R)$ is defined as
$$l(R):=\min\{(-K_X)\cdot C \ | \ C \textrm{ is a rational curve such
that } [C]\in R\}.$$ Note that $0<l(R)\leq n+1$ from Theorem
\ref{cone theorem} and the definition of an extremal rational
curve. Let $\mathcal{E}$ be an ample vector bundle of rank $r\geq
2$ on $X$ and let $\Omega (X,\mathcal{E})$ be the set of extremal
rays $R$ such that $(K_X+\det\mathcal{E})\cdot R<0$. Then it
follows from Theorem \ref{cone theorem} that the set $\Omega
(X,\mathcal{E})$ is finite. For any extremal ray $R$ in $\Omega
(X,\mathcal{E})$ we define a positive integer
$$\Lambda (X,\mathcal{E},R):=(-K_X-\det\mathcal{E})\cdot C,$$
where $C$ is an extremal rational curve such that $-K_X\cdot
C=l(R)$.

\section{Proof of Theorem \ref{thm}.}

Suppose that $K_X+\det\mathcal{E}$ is not nef. By Theorem
\ref{cone theorem} we can find an extremal ray $R$ with
$(K_X+\det\mathcal{E})\cdot R<0$, and so $\Omega
(X,\mathcal{E})\neq\emptyset$. Since the set $\Omega
(X,\mathcal{E})$ is finite, define the positive integer
$$\Lambda (X, \mathcal{E}):=\max \{\Lambda (X, \mathcal{E},R) \ | \
R\in \Omega (X, \mathcal{E}) \}.$$ Therefore we have only the
following three possibilities:

\medskip

\begin{enumerate}
\item[$(1)$] $\Lambda (X, \mathcal{E})\geq 3$; \quad $(2)$
$\Lambda (X, \mathcal{E})=2$; \quad $(3)$ $\Lambda (X,
\mathcal{E})=1$.
\end{enumerate}

\medskip

\noindent From now on, we proceed with a case-by-case analysis.

\medskip

\noindent\textit{Case $(1)$}. Note that, by definition of $\Lambda
(X,\mathcal{E})$, there exists an extremal rational curve $C\in R$
such that $-(K_X+\det\mathcal{E})\cdot C\geq 3$ and
$l(R)=-K_X\cdot C$. So, by Lemma \ref{lemma} we get
$$-K_X\cdot C\geq\det\mathcal{E}\cdot C+3\geq (n-3)+3=n,$$ i.e. $l(R)\geq n$. Let
$f_R:X\to Y$ be the extremal ray contraction associated to $R$.
From \cite[lemma 6.3.12]{BS}, we deduce that $f_R$ is of fibre
type with $\mathrm{dim } Y\leq 1$. Let $F$ be a smooth general
fibre of $f_R$. Thus we have the following two cases:

\medskip

\begin{enumerate}
\item[(a)] $\mathrm{dim } F=n-1, \mathrm{dim } Y=1$; \ \ \ (b)
$\mathrm{dim } F=n, \mathrm{dim } Y=0$.
\end{enumerate}

\medskip

\noindent In Case (a), by \cite[lemma 6.3.12(2)]{BS} we know that
$Y$ is a smooth irreducible curve and $\rho(F)=1$. Note that
$\mathcal{E}|_F$ is an ample vector bundle on $F$ of rank
$\mathrm{dim } F-2$.

\bigskip

\noindent\textit{Claim.} $K_F+\det\mathcal{E}|_F$ is not nef and
$\Lambda (F,\mathcal{E}|_F)\geq 3$ for a general fibre $F$ of
$f_R:X\to Y$.

\medskip

\noindent Since $\rho(F)=1$, we have $\overline{NE}(F)\cong\langle
C'\rangle$ with $[C']\in R=\mathbb{R}_+[C]$, i.e. $C'\equiv\delta
C$ for some $\delta >0$. Thus for any effective curve
$\gamma\subset F$, we know that $\gamma\equiv\delta' C'$ for some
$\delta'>0$ and so
$$(K_F+\det\mathcal{E}|_F)\cdot\gamma
=(K_X+\det\mathcal{E})_F\cdot (\delta'
C')=\delta'(K_X+\det\mathcal{E})\cdot (\delta
C)=\delta\delta'(K_X+\det\mathcal{E})\cdot C<0.$$ This shows that
$K_F+\det\mathcal{E}|_F$ cannot be nef and that $-K_F$ is ample on
$F$. Then by Theorem \ref{cone theorem}, we see that
$\overline{NE}(F)\cong\mathbb{R}_+[C']$, where $C'$ is an extremal
rational curve on $F$. Moreover, since
$$\delta(-K_X\cdot C)=-K_X\cdot C'\geq l(R)=-K_X\cdot C>0,$$ we
deduce that $\delta\geq 1$. Thus we get for $C'\subset F$
$$(-K_F-\det\mathcal{E}|_F)\cdot C'=(-K_X-\det\mathcal{E})\cdot
C'= (-K_X-\det\mathcal{E})\cdot\delta C\geq 3\delta\geq 3,$$ and
this gives $\Lambda(F,\mathcal{E}|_F)\geq 3$ for a general fibre
$F$ of $f_R:X\to Y$. \hfill{Q.E.D.}

\bigskip

\noindent By the \textit{Claim} and \cite[Proposition
1.1']{Zhang}, we obtain that $(F,\mathcal{E}|_F)$ is one of the
following pairs:

\smallskip

\begin{enumerate}
\item[$(a_1)$] $(\mathbb{P}^{n-1},
\mathcal{O}_{\mathbb{P}^{n-1}}(1)^{\oplus n-3})$; \item[$(a_2)$]
$(\mathbb{P}^{n-1},
\mathcal{O}_{\mathbb{P}^{n-1}}(2)\oplus\mathcal{O}_{\mathbb{P}^{n-1}}(1)^{\oplus
n-4})$; \item[$(a_3)$] $(\mathbb{Q}^{n-1},
\mathcal{O}_{\mathbb{Q}^{n-1}}(1)^{\oplus n-3})$.
\end{enumerate}

\smallskip

\noindent In Case $(a_2)$, we get
$-(K_F+\det\mathcal{E}|_F)=\mathcal{O}_{\mathbb{P}^{n-1}}(2)$.
Take a line $l\subset F$. Since $l\in R=\mathbb{R}_+[C]$ and
$l(R)=-K_X\cdot C$, by arguing as in the proof of \textit{Claim
1}, we deduce that $l\equiv\rho C$ for some $\rho\geq 1$, but this
gives the following numerical contradiction
$$2=-(K_F+\det\mathcal{E}|_F)\cdot l=-(K_X+\det\mathcal{E})\cdot\rho
C\geq 3\rho\geq 3.$$ By a similar argument, we see that also Case
$(a_3)$ cannot occur. Consider now Case $(a_1)$ and let
$c:=\det\mathcal{E}\cdot C$. As in Case $(a_2)$, taking a line
$l\subset F$, we deduce that $3=(-K_F-\det\mathcal{E}|_F)\cdot
l\geq -(K_X+\det\mathcal{E})\cdot C,$ i.e.
$(K_X+\det\mathcal{E})\cdot C=-3$. Therefore
$$[cK_X+(c+3)\det\mathcal{E}]\cdot C=c[K_X+\det\mathcal{E}]\cdot
C+3c=0$$ and from \cite[theorem 4.3.1]{BS} it follows that
$$cK_X+(c+3)\det\mathcal{E}\simeq f_R^*D,$$ for some $D\in\mathrm{Pic}(Y)$.
Then we have
\begin{equation}
cK_X+(c+3)(\det\mathcal{E}+(n-3)f_R^*A)\simeq f_R^*H \label{a1}
\end{equation}
for a suitable very ample line bundle $A$ on $Y$ such that
$H:=D+(c+3)(n-3)A$ is ample on $Y$. Restricting (\ref{a1}) to a
general fibre $F$ of $f_R$, we get
$$\mathcal{O}_F\simeq cK_F+(c+3)\det\mathcal{E}|_F=\mathcal{O}_F(-nc+(c+3)(n-3)).$$
Hence $c=n-3$ and (\ref{a1}) becomes
\begin{equation}
(n-3)K_X+n(\det\mathcal{E}+(n-3)f_R^*A)\simeq f_R^*H. \label{a2}
\end{equation}
Put $\mathcal{E}':=\mathcal{E}\otimes f_R^*A$ and note that
$\mathcal{E}'$ is ample and such that
$\det\mathcal{E}'=\det\mathcal{E}+(n-3)f_R^*A$. Thus we get
$$(n-3)K_X+n(\det\mathcal{E}')\simeq f_R^*H.$$

\noindent Moreover, since
$$\left[ K_X+\left(\frac{n}{n-3}\right)\det\mathcal{E}'\right]\cdot
l=\left[
K_F+\left(\frac{n}{n-3}\right)\det\mathcal{E}|_F\right]\cdot
l=0,$$ for any line $l\subset F\cong\mathbb{P}^{n-1}$, we see that
the nef value $\tau(X,\mathcal{E}')$ of the pair
$(X,\mathcal{E}')$ is such that
$\tau(X,\mathcal{E}')\geq\frac{n}{n-3}=\frac{n}{\mathrm{rk
}\mathcal{E}'}>1$. Thus by \cite[proposition 4(3)]{AN} we have
$X\cong\mathbb{P}_Y(\mathcal{V})$ for a suitable rank $n$ vector
bundle $\mathcal{V}$ on $Y$ and from \cite[theorem 1.3(5)]{Ohno}
we conclude that $\mathcal{E}'\cong\xi_{\mathcal{V}}\otimes
f_R^{*}\mathcal{G}$, where $\mathcal{G}$ is an ample vector bundle
of rank $n-3$ on $Y$.

\medskip

Finally, in Case (b), since $X$ is a Fano $n$-fold with
$\rho(X)=1$, for any effective rational curve $\Gamma$ on $X$ we
see that $\Gamma\in R$ and $-K_X\cdot\Gamma\geq l(R)\geq n$. Thus
by \cite{Miyaoka}, \cite{Miyaoka1} and \cite{CMS} we can conclude
that $X$ is either $\mathbb{P}^n$ or $\mathbb{Q}^n$. In the former
case, since $\Lambda(X,\mathcal{E})\geq 3$, by Lemma \ref{lemma}
we have
$$n-3\leq c_1(\mathcal{E})\cdot l\leq -K_{\mathbb{P}^n}\cdot
l-3=n-2,$$ where $l$ is a line in $\mathbb{P}^n$. Therefore by
\cite[theorems 3.2.1 and 3.2.3]{OSS}, $\mathcal{E}$ is isomorphic
to either $\mathcal{O}(1)^{\oplus n-3}$, or
$\mathcal{O}(2)\oplus\mathcal{O}(1)^{\oplus n-4}$. In the latter
case, since $X\cong\mathbb{Q}^n$ and $\Lambda(X,\mathcal{E})\geq
3$, we see that
$$n-3\leq c_1(\mathcal{E})\cdot l\leq -K_{\mathbb{Q}^n}\cdot
l-3=n-3,$$ i.e. $c_1(\mathcal{E})\cdot l=n-3$, and then
$\mathcal{E}\cong\mathcal{O}(1)^{\oplus n-3}$ by \cite[lemma
3.6.1]{W}.

\medskip

Summing up the above discussion, we obtain the following

\begin{Prop}\label{prop 1}
If $\mathcal{E}$ is an ample vector bundle of rank $n-3\geq 2$
such that $K_X+\det\mathcal{E}$ is not nef and
$\Lambda(X,\mathcal{E})\geq 3$, then the pair $(X,\mathcal{E})$ is
one of the following:
\begin{enumerate}
\item[$(1')$] $(\mathbb{P}^n,
\mathcal{O}_{\mathbb{P}^n}(1)^{\oplus n-3})$; \ $(2')$
$(\mathbb{P}^n,
\mathcal{O}_{\mathbb{P}^n}(2)\oplus\mathcal{O}_{\mathbb{P}^n}(1)^{\oplus
n-4})$; \ $(3')$ $(\mathbb{Q}^n,
\mathcal{O}_{\mathbb{Q}^n}(1)^{\oplus n-3})$; \item[$(4')$] there
is a vector bundle $\mathcal{V}$ on a smooth curve $C$ such that
$X\cong\mathbb{P}_C(\mathcal{V})$ and
$\mathcal{E}\cong\xi_{\mathcal{V}}\otimes\pi^*\mathcal{G}'$, where
$\pi:X\to C$ is the projection map, $\xi_{\mathcal{V}}$ is the
tautological line bundle of $\mathbb{P}(\mathcal{V})$ and
$\mathcal{G}'$ is a suitable vector bundle of rank $n-3$ on $C$;
in particular,
$\mathcal{E}|_F\cong\mathcal{O}_{\mathbb{P}^{n-1}}(1)^{\oplus
n-3}$ for any fibre $F\cong\mathbb{P}^{n-1}$ of the map $\pi$.
\end{enumerate}
\end{Prop}

\smallskip

\noindent\textit{Case $(2)$}. In this situation, we know that
there exists an extremal rational curve $C\in R$ such that
$-(K_X+\det\mathcal{E})\cdot C=2$ and $-K_X\cdot C=l(R)$.
Moreover, note that
$$
\left[K_X+\left(1+\frac{2}{\det\mathcal{E}\cdot
C}\right)\det\mathcal{E}\right]\cdot C=0,
$$
i.e. the nef value $\tau(X,\mathcal{E})$ of the pair
$(X,\det\mathcal{E})$ is such that
\begin{equation}
\tau(X,\mathcal{E})\geq 1+\frac{2}{\det\mathcal{E}\cdot C}.
\label{a3}
\end{equation}
Since $-K_X\cdot C\leq n+1$, by Lemma \ref{lemma} we get
$\det\mathcal{E}\cdot C=l(R)-2$ and
$$n-3\leq\det\mathcal{E}\cdot C=-K_X\cdot C-2\leq n-1.$$ So we have
only the following three possibilities:

\smallskip

\begin{enumerate}
\item[(I)] $\det\mathcal{E}\cdot C=n-1,
\tau(X,\mathcal{E})\geq\frac{n+1}{n-1}=\frac{l(R)}{n-1}$;
\item[(II)] $\det\mathcal{E}\cdot C=n-2,
\tau(X,\mathcal{E})\geq\frac{n}{n-2}=\frac{l(R)}{n-2}$;
\item[(III)] $\det\mathcal{E}\cdot C=n-3,
\tau(X,\mathcal{E})\geq\frac{n-1}{n-3}=\frac{l(R)}{n-3}$.
\end{enumerate}

\smallskip

\noindent In Case (I), we see that $-K_X\cdot C=n+1$, i.e.
$l(R)=n+1$. From the Ionescu-Wi\'{s}niewski inequality,
\cite[lemma 6.3.12]{BS} and \cite{CMS}, it follows that
$X\cong\mathbb{P}^n$. Since by (\ref{a3}) we have
$$\tau(X,\mathcal{E})\geq\frac{n+1}{n-1}\geq\frac{n-2}{n-3}=\frac{n-2}{\mathrm{rk}\ \mathcal{E}}>1$$
and $\det\mathcal{E}\cdot C-(n-3)=2$, we obtain by
\cite[proposition 6(a)(1)]{AN} and \cite{Ohno} that $\mathcal{E}$
is a decomposable vector bundle with
$c_1(\mathcal{E})=\mathcal{O}_{\mathbb{P}^n}(n-1)$. Since
$\mathcal{E}$ is ample, this gives
$\mathcal{E}\cong\mathcal{O}(3)\oplus\mathcal{O}(1)^{\oplus n-4}$
or $\mathcal{O}(2)^{\oplus 2}\oplus\mathcal{O}(1)^{\oplus n-5}$.

\medskip

\noindent In Case (II), note that
$$\tau(X,\mathcal{E})\geq\frac{n}{n-2}>\frac{n-2}{n-3}=\frac{n-2}{\mathrm{rk}\ \mathcal{E}}>1$$ and
$\det\mathcal{E}\cdot C-(n-3)=1$. Moreover, $l(R)=-K_X\cdot C=n$
and from \cite[lemma 6.3.12]{BS} we have only the following two
cases:

\smallskip

\begin{enumerate}
\item[$(A)$] $X$ is a Fano $n$-fold with $\rho(X)=1$; \item[$(B)$]
$\rho(X)=2$ and the contraction $f_R:X\to Y$ associated to the
extremal ray $R=\mathbb{R}_+[C]$ is a morphism onto a smooth curve
$Y$, whose general fibre $F$ is a smooth Fano $(n-1)$-fold with
$\rho(F)=1$.
\end{enumerate}

\smallskip

\noindent In Case $(A)$, by \cite{Miyaoka} and \cite{Miyaoka1} we
obtain that $X\cong\mathbb{Q}^n\subset\mathbb{P}^{n+1}$ with
$c_1(\mathcal{E})=\mathcal{O}_{\mathbb{Q}^n}(n-2)$. Thus from
\cite[theorem 1.3(15)]{Ohno} we get
$\mathcal{E}\cong\mathcal{O}_{\mathbb{Q}^n}(2)\oplus\mathcal{O}_{\mathbb{Q}^n}(1)^{\oplus
n-4}$. In Case $(B)$, we conclude by \cite[theorem 1.3(10)]{Ohno}.

\medskip

\noindent Finally, in Case (III) by \cite[theorem 1]{AN}
we have
$$\tau(X,\mathcal{E})=\frac{n-1}{n-3}=\frac{l(R)}{\mathrm{rk}\ \mathcal{E}}.$$
Since $l(R)=-K_X\cdot C=n-1$ and $n=r+3\geq 5$, we can conclude by
\cite[propositions $3,4$ and $5$]{AN} and \cite[theorem 1.3(6),
(8), (11) and (12)(a)]{Ohno}. Moreover, if $\rho(X)=1$ then from
\cite[theorem 1(2)]{AN} it follows that
$\mathcal{E}\cong\mathcal{O}_X(1)^{\oplus n-3}$, which leads to
Case $(5)(\alpha)$ of Theorem \ref{thm}.

\medskip

The above discussion can be summed up in the following

\begin{Prop}\label{prop 2}
If $\mathcal{E}$ is an ample vector bundle of rank $n-3\geq 2$
such that $K_X+\det\mathcal{E}$ is not nef and
$\Lambda(X,\mathcal{E})=2$, then the pair $(X,\mathcal{E})$ is one
of the following:
\begin{enumerate}
\item[$(1'')$] $(\mathbb{P}^n,
\mathcal{O}_{\mathbb{P}^n}(3)\oplus\mathcal{O}_{\mathbb{P}^n}(1)^{\oplus
n-4})$; \item[$(2'')$] $(\mathbb{P}^n,
\mathcal{O}_{\mathbb{P}^n}(2)^{\oplus
2}\oplus\mathcal{O}_{\mathbb{P}^n}(1)^{\oplus n-5})$;
\item[$(3'')$] $(\mathbb{Q}^n,
\mathcal{O}_{\mathbb{Q}^n}(2)\oplus\mathcal{O}_{\mathbb{Q}^n}(1)^{\oplus
n-4})$; \item[$(4'')$] there exist a vector bundle $\mathcal{V}$
on a smooth curve $C$ such that $X\cong\mathbb{P}_C(\mathcal{V})$
and an exact sequence
$$0\to\pi^*\mathcal{L}\otimes\xi_{\mathcal{V}}^{\otimes
2}\to\mathcal{E}\to\pi^*\mathcal{G}\otimes\xi_{\mathcal{V}}\to 0$$
for some line bundle $\mathcal{L}$ on $C$ and some vector bundle
$\mathcal{G}$ of rank $n-4\geq 1$ on $C$, where
$\pi:\mathbb{P}(\mathcal{V})\to C$ is the projection map and
$\xi_{\mathcal{V}}$ is the tautological line bundle on
$\mathbb{P}_C(\mathcal{V})$; in particular, we have
$\mathcal{E}|_F\cong\mathcal{O}_{\mathbb{P}^{n-1}}(2)\oplus\mathcal{O}_{\mathbb{P}^{n-1}}(1)^{\oplus
n-4}$ for any fibre $F\cong\mathbb{P}^{n-1}$ of the map $\pi$;
\item[$(5'')$] $X$ is a Del Pezzo $n$-fold with $\rho(X)=1$ and
$\mathcal{E}\cong L^{\oplus n-3}$, where $L$ is the ample
generator of $\mathrm{Pic }(X)$; \item[$(6'')$] $X$ is a quadric
fibration, $q:X\to C$, of the relative Picard number one over a
smooth curve $C$ (i.e. $X$ is a section of a divisor of relative
degree two in a projective space $\mathbb{P}_C(\mathcal{F})$,
where $\mathcal{F}$ is a vector bundle of rank $n+1$ on $C$) and
there exist a $q$-ample line bundle $\mathcal{O}_X(1)$ on $X$ and
an ample vector bundle $\mathcal{G}$ of rank $n-3$ on $C$ such
that $\mathcal{E}\cong\mathcal{O}_X(1)\otimes q^*\mathcal{G}$,
where $\mathcal{O}_X(1)|_F\simeq\mathcal{O}_{\mathbb{Q}^{n-1}}(1)$
for any fibre $F\cong\mathbb{Q}^{n-1}\subset\mathbb{P}^{n}$ of the
map $q$; \item[$(7'')$] there exists a
$\mathbb{P}^{n-2}$-fibration $\pi:X\to S$, locally trivial in the
complex (or \'{e}tale) topology, over a smooth surface $S$ such
that $\mathcal{E}|_F\cong\mathcal{O}_{\mathbb{P}^{n-2}}(1)^{\oplus
n-3}$ for any closed fiber $F\cong\mathbb{P}^{n-2}$ of the map
$\pi$; \item[$(8'')$] there exists a smooth variety $X'$ and a
morphism $\varphi:X\to X'$ expressing $X$ as the blowing-up of
$X'$ at a finite set of points $B$ and an ample vector bundle
$\mathcal{E}'$ on $X'$ such that
$\mathcal{E}\otimes([\varphi^{-1}(B)])\cong\varphi^*\mathcal{E}'$
and $K_{X'}+\tau'\det\mathcal{E}'$ is nef, where $\tau' <
\frac{n-1}{n-3}$ is the nefvalue of the pair
$(X',\det\mathcal{E}')$. Moreover,
$\mathcal{E}|_E\cong\mathcal{O}_{\mathbb{P}^{n-1}}(1)^{\oplus
n-3}$, where $E$ is any irreducible component of the exceptional
locus of the map $\varphi$.
\end{enumerate}
\end{Prop}

\smallskip

\noindent\textit{Case $(3)$}. By arguing as in \cite[p.77]{Maeda}
and \cite[proposition 3.5]{W}, note that in fact we get $L\cdot
C=1$, where $L:=(m-1)K_X+m\det\mathcal{E}$ is ample on $X$ and
$m\geq n-3\geq 2$ with
$$m:=\min\{ (\det\mathcal{E})\cdot C \ | \ C \textrm{ extremal rational curve s.t. }
(K_X+\det\mathcal{E})\cdot C=-1\}.$$ So, define
$\mathcal{V}:=\mathcal{E}\oplus L$. Then $\mathcal{V}$ is an ample
vector bundle of rank $n-2$ on $X$ such that
$$(K_X+\det\mathcal{V})\cdot C=(K_X+\det\mathcal{E})\cdot C+L\cdot
C=-1+1=0.$$ Therefore $K_X+\det\mathcal{V}$ cannot be ample on
$X$. Note also that $K_X+\det\mathcal{V}$ must be nef. Otherwise,
by Theorem \ref{cone theorem} we could find an extremal ray $R$ on
$X$ such that
$$-1\geq (K_X+\det\mathcal{V})\cdot R\geq
(K_X+\det\mathcal{E})\cdot R+1,$$ i.e. $(K_X+\det\mathcal{E})\cdot
R\leq -2$, but this would contradict $\Lambda(X,\mathcal{E})=1$.

Thus assume that $K_X+\det\mathcal{V}$ is nef but not big. From
\cite[theorem 5.1(2)]{AM} we know that there exists a morphism
$\Phi:X\to W$ onto a normal variety $W$ supported by (a large
multiple of) $K_X+\det\mathcal{V}$ and $\mathrm{dim } W\leq 3$.
Moreover, let $F$ be a general fibre of $\Phi$. We have the
following according to $s:=\mathrm{dim } W$:

\smallskip

\begin{enumerate}
\item[($i_1$)] if $s=0$ then $X$ is a Fano $n$-fold with
$K_X+\det\mathcal{V}=\mathcal{O}_X$. If $n\geq 6$, then $b_2(X)=1$
except if $X\cong\mathbb{P}^3\times\mathbb{P}^3$ and
$\mathcal{E}\cong\mathcal{O}(1,1)^{\oplus 4}$; \item[($i_2$)] if
$s=1$ then $W$ is a smooth curve and $\Phi$ is a flat
(equidimensional) map. If $n\geq 6$ then $\Phi$ is an elementary
contraction. Moreover, by \cite{Occhetta} and \cite{PSW},
$(F,\mathcal{V}|_F)$ is one of the following pairs:
\begin{enumerate}
\item[$(a)$] $(\mathbb{P}^{n-1},
\mathcal{O}_{\mathbb{P}^{n-1}}(2)^{\oplus
2}\oplus\mathcal{O}_{\mathbb{P}^{n-1}}(1)^{\oplus n-4})$;
\item[$(b)$] $(\mathbb{P}^{n-1},
\mathcal{O}_{\mathbb{P}^{n-1}}(3)\oplus\mathcal{O}_{\mathbb{P}^{n-1}}(1)^{\oplus
n-3})$; \item[$(c)$] $(\mathbb{Q}^{n-1},
\mathcal{O}_{\mathbb{Q}^{n-1}}(2)\oplus\mathcal{O}_{\mathbb{Q}^{n-1}}(1)^{\oplus
n-3})$; \item[$(d)$] $(\mathbb{Q}^{4},
\mathcal{S}(2)\oplus\mathcal{O}_{\mathbb{Q}^{n-1}}(1))$, where
$\mathcal{S}$ is the spinor bundle on
$\mathbb{Q}^4\subset\mathbb{P}^5$; \item[$(e)$] $F$ is a Del Pezzo
$(n-1)$-fold with $\rho(F)=1$, $-K_F\simeq (n-2)H_F$,
Pic$(F)\cong\mathbb{Z}[H_F]$ and $\mathcal{V}|_F\cong H_F^{\oplus
n-2}$; \item[$(f)$]
$(\mathbb{P}^2\times\mathbb{P}^2,\mathcal{O}(1,1)^{\oplus 3})$.
\end{enumerate}
If the general fibre is $\mathbb{P}^{n-1}$ then $X$ is a classical
scroll, while if the general fibre is
$\mathbb{Q}^{n-1}\subset\mathbb{P}^n$ then $X$ is a quadric
bundle; \item[($i_3$)] if $s=2$ then $W$ is a smooth surface,
$\Phi$ is a flat map and by \cite{F1} and \cite{P1},
$(F,\mathcal{V}|_F)$ is one of the following pairs:
\begin{enumerate}
\item[$(\alpha)$] $(\mathbb{P}^{n-2},
\mathcal{O}_{\mathbb{P}^{n-2}}(2)\oplus\mathcal{O}_{\mathbb{P}^{n-2}}(1)^{\oplus
n-3})$; \item[$(\beta)$] $(\mathbb{P}^{n-2},
T_{\mathbb{P}^{n-2}})$, where $T_{\mathbb{P}^{n-2}}$ is the
tangent bundle on $\mathbb{P}^{n-2}$; \item[$(\gamma)$]
$(\mathbb{Q}^{n-2}, \mathcal{O}_{\mathbb{Q}^{n-2}}(1)^{\oplus
n-2})$.
\end{enumerate}
In particular, if the general fibre is $\mathbb{P}^{n-2}$ then all
the fibres are $\mathbb{P}^{n-2}$; \item[($i_4$)] if $s=3$ then
$W$ is a smooth $3$-fold and all the fibres of $\Phi$ are
isomorphic to $\mathbb{P}^{n-3}$.
\end{enumerate}

\smallskip

\noindent Furthermore, since $(K_X+\det\mathcal{V})\cdot C=0$,
note that $L_F\cdot C=L\cdot C=1$, $C$ being an extremal rational
curve of $X$ which belongs into a fibre of $\Phi$. Moreover, from
$(K_X+\det\mathcal{E})\cdot C=-1$, we deduce that
$$-K_X\cdot C=\det\mathcal{E}\cdot C+1\geq m+1=(m+1)L\cdot C,$$
and then the nef value $\tau(X,L)$ of the polarized pair $(X,L)$
is greater than or equal to $m+1\geq n-2$.

\medskip

\noindent \textit{Case} ($i_1$). Observe that $X$ is a Fano
$n$-fold with $-K_X\simeq\eta L$ and
$$\eta=\tau(X,L)\geq m+1\geq n-2\geq 3.$$ If $\rho(X)=1$, then
by \cite[chapter 7]{BS} we know that $X$ is either $\mathbb{P}^n,
\mathbb{Q}^n$, a Del Pezzo $n$-fold with $-K_X\simeq (n-1)L$, or a
Mukai $n$-fold with $-K_X\simeq (n-2)L$. Since
$\det\mathcal{E}\cdot C=-K_X\cdot C-1\geq\eta-1$, note that $m\geq
\eta-1$, that is, $\eta\leq m+1$. Thus $\eta =\tau(X,L)=m+1$ and
so we get
$$L=(\eta-2)K_X+(\eta-1)\det\mathcal{E}=-\eta (\eta-2)L+(\eta-1)\det\mathcal{E},$$
i.e. $\det\mathcal{E}=(\eta-1)L$, where
Pic$(X)\cong\mathbb{Z}[L]$. Let us proceed with a case--by--case
analysis.

\smallskip

First of all, assume that $X\cong\mathbb{P}^n$. Note that
$L=\mathcal{O}_{\mathbb{P}^n}(1)$ and
$c_1(\mathcal{E})=\mathcal{O}_{\mathbb{P}^n}(n)$. Without loss of
generality, we can write
$\mathcal{E}\cong V_1\oplus V_2$, where $V_1$ is an indecomposable
vector bundle on $\mathbb{P}^n$. Put $r_i=\mathrm{rk} V_i$ and
note that $r_1+r_2=n-3$ with $1\leq r_1\leq n-3$. Furthermore, if
$\deg V_i\leq r_i+1$, then $V_i$ is an uniform vector bundle on
$\mathbb{P}^n$ of rank $r_i\leq n-3$, i.e. $V_i$ is a sum of line
bundles on $\mathbb{P}^n$. Moreover, since $V_1$ is an ample
vector bundle on $\mathbb{P}^n$ such that $\deg V_1\leq
\deg\mathcal{E}=n$, if $r_1=2$ then by \cite[theorem (9.1)]{APW}
we see that $V_1$ splits. Thus we can assume that $r_1\geq 3$.
Then we get the following three possibilities:

\medskip

(a) $r_1=n-3\geq 3$ and $\deg V_1=n$;

\smallskip

(b) $3\leq r_1\leq n-4$, $\deg V_1=r_1+2$ and $\deg V_2=r_2+1$;

\smallskip

(c) $3\leq r_1\leq n-4$, $\deg V_1=r_1+3$ and $\deg V_2=r_2$.

\bigskip

\noindent By Lemma \ref{lemma 1}, we have $V_1(-1)$ is a nef
vector bundle on $\mathbb{P}^n$ either in Case (b), or in Case (c)
for $3\leq r_1\leq n-6$. Note that Case (b) does not occur.
Indeed, since $V_1(-1)$ is nef and $c_1(V_1(-1))=2$, by
\cite[Lemma 6 (1)]{PSW1} we know that $V_1(-1)$ is decomposable,
but this leads to a contradiction. In Case (c), since $V_2$ is
ample, we get that $V_2\cong\mathcal{O}_{\mathbb{P}^n}(1)^{\oplus
r_2}$, that is,
$\mathcal{E}\cong\mathcal{O}_{\mathbb{P}^n}(1)^{\oplus r_2}\oplus
V_1$. Finally, assume that $V_1(-1)$ is globally generated on
$\mathbb{P}^n$. Since $n\geq 6$, from Corollary \ref{corollary 3}
we know that $V_1$ splits, a contradiction. This completes Case
$(2)$ of Theorem \ref{thm}.

\smallskip

Suppose now that $X\cong\mathbb{Q}^n\subset\mathbb{P}^{n+1}$. Note
that $L=\mathcal{O}_{\mathbb{Q}^n}(1)$ and
$c_1(\mathcal{E})=\mathcal{O}_{\mathbb{Q}^n}(n-1)$. As above,
write $\mathcal{E}\cong V_1\oplus V_2$, where $V_1$ is an
indecomposable vector bundle on $\mathbb{Q}^n$. Put
$r_i=\mathrm{rk} V_i$ and note that $r_1+r_2=n-3$ with $1\leq
r_1\leq n-3$. Furthermore, if $\deg V_i\leq r_i+1$, then $V_i$ is
an uniform vector bundle on $\mathbb{Q}^n$ of rank $r_i\leq n-3$
and by \cite{KS} we know that $V_i$ is a sum of line bundles on
$\mathbb{Q}^n$. Thus $\deg V_1\geq r_1+2$ and then
$$\deg V_2=n-1-\deg V_1\leq n-3-r_1=r_2.$$
This shows that $\deg V_2=r_2$ and $\deg V_1=r_1+2$, that is,
$\mathcal{E}\cong\mathcal{O}_{\mathbb{Q}^n}(1)^{\oplus
n-3-r_1}\oplus V_1$. Since $V_1$ is an ample vector bundle on
$\mathbb{Q}^n$ such that $\deg V_1\leq\deg\mathcal{E}=n-1$, we
deduce that $V_1$ is a Fano bundle on $\mathbb{Q}^n$. Thus if
$r_1=2$, then by \cite[theorem (2.4)(2)]{APW} we see that either
$V_1$ splits, or $n=5$ and $V_1=\mathcal{C}(a)$, where
$\mathcal{C}$ is a Cayley bundle on $\mathbb{Q}^5$ such that
$c_1(\mathcal{C})=-1$. Since $4=c_1(V_1)=-1+2a$, we see that $a$
is not an integer, but this is not possible. Hence $V_1$ splits
for $r_1=2$. On the other hand, if $r_1\geq 3$ then by restricting
$V_1$ to any line on $\mathbb{Q}^n$, we have the only possible
splitting types are $(3,1,1,...,1)$ or $(2,2,1,...,1)$. So either
$V_1$ is a uniform vector bundle of splitting type $(3,1,1,...,1)$
and by \cite{KS} $V_1$ is a direct sum of line bundles on
$\mathbb{Q}^n$, or $V_1$ has the generic splitting type
$(2,2,1,...,1)$. Finally, by Lemma \ref{lemma 1} and Corollary
\ref{corollary 5}, we know that $V_1(-1)$ is nef for $3\leq
r_1\leq n-6$, and it is not globally generated for $n\geq 7$. This
gives Case $(4)$ of Theorem \ref{thm}.

\smallskip

Assume that $X$ is a Fano $n$-fold of index $n-1$. Since
$c_1(\mathcal{E})=\mathcal{O}_X(n-2)$, we see that for every line
$l$ of $(X,H)$, where $H$ is the ample generator of
$\mathrm{Pic}(X)$, it happens that $\mathcal{E}|_l\cong
H|_l^{\oplus n-4}\oplus 2H|_l$, i.e. $\mathcal{E}$ is a uniform
vector bundle on $X$ with splitting type $(2,1,1,...,1)$. This
leads to Case $(5)(\beta)$ of Theorem \ref{thm}.

\smallskip

Finally, suppose that $X$ is a Fano $n$-fold of index $n-2$. From
\cite[theorem 2.10]{Mella}, it follows by hyperplane sections that
there exists at least a line $l$ of $(X,L)$, where
$\mathrm{Pic}(X)=\mathbb{Z}[L]$. Recall that
$\mathcal{V}=\mathcal{E}\oplus L$. Since $\rho(X)=1$ and
$$\tau(X,\mathcal{V})\geq 1=\frac{n-2}{\mathrm{rk}\mathcal{V}}=\frac{-K_X\cdot l}{\mathrm{rk}\mathcal{V}},$$
by \cite[proposition 1.2]{AW} and \cite[theorem 1(2)]{AN}, we
obtain that $\tau(X,\mathcal{V})=\frac{-K_X\cdot
l}{\mathrm{rk}\mathcal{V}}$ and so $\mathcal{V}\cong L^{\oplus
n-2}$, i.e. $\mathcal{E}\cong L^{\oplus n-3}$. This gives Case
$(6)$ in Theorem \ref{thm}.

\medskip

Assume now that $\rho(X)\geq 2$. Since $\eta\geq
n-2\geq\frac{n+1}{2}$, from \cite{Occhetta1}, \cite{W1} and
\cite{W2}, it follows that $\eta\leq\frac{n+2}{2}$, i.e. $5\leq
n\leq 6$, and that $(X,\mathcal{V},L)$ is one of the following
triplets:
\begin{itemize}
\item[-]
$(\mathbb{P}^3\times\mathbb{P}^3,\mathcal{O}_{\mathbb{P}^3\times\mathbb{P}^3}(1,1)^{\oplus
4},\mathcal{O}_{\mathbb{P}^3\times\mathbb{P}^3}(1,1))$; \item[-]
$(\mathbb{P}^2\times\mathbb{Q}^3,\mathcal{O}_{\mathbb{P}^2\times\mathbb{Q}^3}(1,1)^{\oplus
3}, \mathcal{O}_{\mathbb{P}^2\times\mathbb{Q}^3}(1,1))$; \item[-]
$(\mathbb{P}_{\mathbb{P}^3}(T),[\xi_{T}]^{\oplus 3}, \xi_{T})$,
where $T$ is the tangent bundle of $\mathbb{P}^3$ and $\xi_{T}$ is
the tautological line bundle; \item[-]
$(\mathbb{P}_{\mathbb{P}^3}(\mathcal{O}(2)\oplus\mathcal{O}(1)^{\oplus
2}), [\xi]^{\oplus 3}, \xi)$, where $\xi$ is the tautological line
bundle.
\end{itemize}
Since $\mathcal{V}=\mathcal{E}\oplus L$, we obtain Cases $(7)$ to
$(10)$ of Theorem \ref{thm}.

\medskip

\noindent\textit{Case} ($i_2$). Note that
$L|_F=-(K_X+\det\mathcal{E})|_F=-K_F-\det\mathcal{E}|_F$. Then
$$K_F+mL|_F=K_F+m[-K_F-\det\mathcal{E}|_F]=-[(m-1)K_F+m\det\mathcal{E}|_F]=-L|_F,$$ since
$L:=(m-1)K_X+m\det\mathcal{E}$. This gives
$$L|_F=-K_F-\det\mathcal{E}|_F=(m+1)L|_F-\det\mathcal{E}|_F,$$ i.e.
$\det\mathcal{E}|_F=mL|_F$. Moreover, observe also that the
extremal rational curve $C$ such that $L\cdot C=1$ is contained in
a fibre of $\Phi$, since $(K_X+\det\mathcal{V})\cdot C=0$. Keeping
in mind that $\mathcal{V}=\mathcal{E}\oplus L$ and that
$L=(m-1)K_X+m\det\mathcal{E}$ with $m\geq n-3$, by all the above
remarks we can easily obtain Cases $(11)(c)(d),
(12)\mathit{(ii)}\mathit{(iii)}$, and $(13)$ of Theorem \ref{thm}.

\medskip

\noindent\textit{Case} ($i_3$). Note that Case $(\beta)$ cannot
occur, since $\mathcal{V}|_F=\mathcal{E}|_F\oplus L|_F$ is a
decomposable vector bundle on $\mathbb{P}^{n-2}$ with $n-2\geq 3$.
Moreover, since the extremal rational curve $C$ such that $L\cdot
C=1$ is contained in a fibre of $\Phi$, we obtain immediately
Cases $(14)$ and $(16)$ of Theorem \ref{thm}.

\medskip

\noindent\textit{Case} ($i_4$). Since
$K_F+\det\mathcal{V}|_F=\mathcal{O}_F$, where $\mathcal{V}|_F$ is
an ample vector bundle of rank $\mathrm{dim } F+1$, from \cite{P}
we deduce that $(F,\mathcal{V}|_F)\cong
(\mathbb{P}^{n-3},\mathcal{O}_{\mathbb{P}^{n-3}}(1)^{\oplus n-2})$
and then $L|_F=\mathcal{O}_{\mathbb{P}^{n-3}}(1)$ for any fibre
$F\cong\mathbb{P}^{n-3}$ of $\Phi$. Moreover, as observed in
\textit{Case} ($i_2$), the extremal rational curve $C$ such that
$L\cdot C=1$ is contained in a fibre of $\Phi$. Thus by
\cite[(2.12)]{Fujita} we get immediately Case $(17)$ of Theorem
\ref{thm}.

\medskip

Finally, suppose that $K_X+\det\mathcal{V}$ is nef and big, but
not ample. Then from \cite[theorem 5.1(3)]{AM} we deduce that a
high multiple of $K_X+\det\mathcal{V}$ defines a birational map,
$\varphi:X\to X'$, which contracts an extremal face. Let
$R_i=\mathbb{R}_+[C_i]$ be the extremal rays spanning this face
for some $i$ in a finite set of index. Call $\rho_i:X\to W_i$ the
contraction associated to one of the $R_i$. Then we have each
$\rho_i$ is birational and of divisorial type. If $D_i$ is one of
the exceptional divisors and $Z_i=\rho_i(D_i)$, we know that
$\mathrm{dim } Z_i\leq 1$ and one of the following possibilities
can occur:

\smallskip

\begin{enumerate}
\item[$(j_1)$] $\mathrm{dim } Z_i=0, D_i\cong\mathbb{P}^{n-1}$ and
$([D_i]_{D_i},\mathcal{V}|_{D_i})$ is either $(\alpha)$
$(\mathcal{O}(-2),\mathcal{O}(1)^{\oplus n-2})$, or $(\beta)$
$(\mathcal{O}(-1),\mathcal{O}(1)^{\oplus
n-3}\oplus\mathcal{O}(2))$; \item[$(j_2)$] $\mathrm{dim } Z_i=0,
D_i$ is a (possible singular) quadric
$\mathbb{Q}^{n-1}\subset\mathbb{P}^{n}$ and
$[D_i]_{D_i}=\mathcal{O}(-1)$; moreover,
$\mathcal{V}|_{D_i}\cong\mathcal{O}(1)^{\oplus n-2}$;
\item[$(j_3)$] $\mathrm{dim } Z_i=1$, $W_i$ and $Z_i$ are smooth
projective varieties and $\rho_i$ is the blow-up of $W_i$ along
$Z_i$; moreover, $\mathcal{V}|_{F_i}\cong\mathcal{O}(1)^{\oplus
n-2}$ for any fibre $F_i\cong\mathbb{P}^{n-2}$ of $D_i\to Z_i$.
\end{enumerate}

\smallskip

\noindent Furthermore, we know also that $\varphi$ is a
composition of disjoint extremal contractions as in $(j_1), (j_2)$
and $(j_3)$. Note that the extremal rational curves $C_i$ such
that $R_i=\mathbb{R}_+[C_i]$ and $L\cdot C_i=1$ are contained in
$D_i$, since $(K_X+\det\mathcal{V})\cdot C_i=0$. Moreover
\begin{equation}
\det\mathcal{E}\cdot C_i=-K_X\cdot C_i-1=\left\{
\begin{array}{ll}
n-2 & \ \mathrm{\ in\ Case}\ (j_1)(\alpha) \\
n-3 & \ \mathrm{\ otherwise }
\end{array}
\right.  \notag
\end{equation}
This shows that if some $\rho_i$ is a contraction associated to an
extremal ray as in Case $(j_1)(\alpha)$, then $\varphi :X\to X'$
contracts an extremal face spanned by only extremal rays as in
$(j_1)(\alpha)$. So, keeping in mind that
$\mathcal{V}=\mathcal{E}\oplus L, L\cdot C_i=1$ and
$L=(m-1)K_X+m\det\mathcal{E}$ with $m\geq n-3\geq 2$, we get Cases
$(19)$ and $(20)$ of Theorem \ref{thm}.

Taking into account the above discussion and the Propositions
\ref{prop 1} and \ref{prop 2}, we obtain the statement of Theorem
\ref{thm}. \hfill{$\square$}

\medskip

Finally, under the extra assumption that $\rho(X)=1$ and
$\mathcal{E}(-1):=\mathcal{E}\otimes\mathcal{O}_X(-1)$ is globally
generated on $X$, we have the following

\begin{Cor}\label{corollary thm}
Let $\mathcal{E}$ be an ample vector bundle of rank $r\geq 2$ on a
smooth complex projective variety $X$ of dimension $n$ with
$\rho(X)=1$. Assume that $\mathcal{E}(-1)$ is globally generated
on $X$. If $r=n-3$, then $K_X+\det\mathcal{E}$ is nef except when
$(X,\mathcal{E})$ is one of the following pairs:
\begin{enumerate}
\item[$(a)$] $(\mathbb{P}^n, \oplus_{i=1}^{n-3}
\mathcal{O}_{\mathbb{P}^n}(a_i)$, where all the $a_i$'s are
positive integers such that $\sum_{i=1}^{n-3}a_i\leq n$;
\item[$(b)$] $(\mathbb{Q}^n,
\oplus_{j=1}^{n-3}\mathcal{O}_{\mathbb{Q}^n}(b_j))$, where all the
$b_j$'s are positive integers such that $\sum_{j=1}^{n-3}b_j\leq
n-1$; \item[$(c)$] $X\cong\mathbb{Q}^6$ and $\mathcal{E}$ is an
indecomposable Fano vector bundle on $\mathbb{Q}^6$ of rank $3$
such that $c_1(\mathcal{E})=5$ and its generic splitting type is
$(2,2,1)$; \item[$(d)$] $X$ is a Fano $n$-fold of index $n-1$ with
$\mathrm{Pic}(X)$ generated by an ample line bundle $\mathcal{H}$
and either $(\alpha)$ $\mathcal{E}\cong\mathcal{H}^{\oplus n-3}$,
or $(\beta)$ $\mathcal{E}|_l\cong\mathcal{H}_l^{\oplus
n-4}\oplus\mathcal{H}_l^{\otimes 2}$ for every line $l$ of
$(X,\mathcal{H})$; \item[$(e)$] $X$ is a Fano $n$-fold of index
$n-2$ with $\mathrm{Pic}(X)$ generated by an ample line bundle
$\mathcal{L}$ and $\mathcal{E}\cong\mathcal{L}^{\oplus n-3}$;
\end{enumerate}
\end{Cor}

Let us show that in the above result Case $(c)$ really exists.

\medskip

\noindent\textit{Example}. Let $\mathbb{Q}^6\subset\mathbb{P}^7$
be the $6$-dimensional quadric hypersurface endowed with a spinor
bundle $\mathcal{S}$ of rank $4$. Note that
$c_4(\mathcal{S}^{\nu})=c_4(\mathcal{S})=0$ and that
$\mathcal{S}^{\nu}$ is globally generated on $\mathbb{Q}^6$, where
$\mathcal{S}^{\nu}$ is the dual bundle of $\mathcal{S}$ (see
\cite[(2.8)(ii), (2.9)]{Ot}). Consider now the following exact
sequence
$$0\to\mathcal{O}_{\mathbb{Q}^6}\to\mathcal{S}^{\nu}\to\mathcal{F}\to
0$$ given by a section of $\mathcal{S}^{\nu}$. Then $\mathcal{F}$
is a globally generated vector bundle on $\mathbb{Q}^6$ of rank
$3$ such that $c_1(\mathcal{F})=-c_1(\mathcal{S})=2$
(\cite[(2.9)]{Ot}). Put
$\mathcal{E}:=\mathcal{F}\otimes\mathcal{O}_{\mathbb{Q}^6}(1)$.
This gives an example of an ample indecomposable vector bundle
$\mathcal{E}$ of rank $3$ on $\mathbb{Q}^6$ with
$c_1(\mathcal{E})=5$ and such that
$\mathcal{E}\otimes\mathcal{O}_{\mathbb{Q}^6}(-1)\cong\mathcal{S}^{\nu}/\mathcal{O}_{\mathbb{Q}^6}$
is globally generated on $\mathbb{Q}^6$.

\section{Two applications.}

Within the adjunction theory, let us give here two easy and
immediate consequences of Theorem \ref{thm}.

\subsection{Ample vector bundles and special varieties in adjunction theory.}\label{3.1}

Let $\mathcal{E}$ be an ample vector bundle of rank $r\geq 2$ on
an $n$-fold $X$. Given a smooth projective variety $Z$, the
classification of such varieties $X$ containing $Z$ as an ample
divisor occupies an important position in the theory of polarized
varieties. Moreover, it is well-known that the structure of $Z$
imposes severe restrictions on that of $X$. Inspired by this
philosophy, we generalize some results on ample divisors to ample
vector bundles. In particular, by using a peculiarity of adjoint
bundles, Lanteri and Maeda \cite{LM} investigated when $K_Z+tH_Z$
is not nef for $t\geq\dim Z-2$, where $H$ is an ample line bundle
on $X$ and $H_Z:=H|_{Z}$. So, following here this idea, we deal
with the next case $t=\dim Z-3$. More precisely, we give a
detailed description of triplets $(X,\mathcal{E},H)$ as above
under the assumption that $K_Z+(\dim Z-3)H_Z$ is not nef for $\dim
Z\geq 3$.

\smallskip

In particular, for $\dim Z=3$, we obtain the following

\begin{Prop}\label{dim Z=3}
Let $X$ be an $n$-fold. Let $\mathcal{E}$ be an ample vector
bundle of rank $r\geq 2$ on $X$ such that there exists a global
section $s\in\Gamma(\mathcal{E})$ whose zero locus $Z=(s)_0\subset
X$ is a smooth subvariety of dimension $n-r=3$. Then $K_Z$ is nef
if and only if the pair $(X,\mathcal{E})$ is not as in
\em{Theorem} \ref{thm} with $n\geq 5$.
\end{Prop}

\noindent\textbf{Proof.} To prove the `if' part, assume that $K_Z$
is not nef. Recalling that $K_Z=(K_X+\det\mathcal{E})_Z$, we
deduce that $K_X+\det\mathcal{E}$ is not nef. Therefore we have
$(X,\mathcal{E})$ is as in Theorem \ref{thm}.

\smallskip

Conversely, to prove the `only if' part of the statement, assume
that the pairs $(X,\mathcal{E})$ are as in Theorem \ref{thm}.
Therefore, if $(X,\mathcal{E})$ is as in cases $(1)$ to $(10)$, it
is easy to see that $-K_Z$ is ample, and so $K_Z$ cannot be nef.

\smallskip

Suppose that $(X,\mathcal{E})$ is as in case $(11)$. Let $\pi$ be
the projection map $X\to C$ and let $s_F$ denote the restriction
of the global section $s$ to a general fiber $F$ of $\pi$. Then
$s_F\in\Gamma(\mathcal{E}|_F)$. Putting $D:=(s_F)_0=Z\cap F$, we
have $\dim D\geq 2$. On the other hand, since $Z$ is irreducible,
we have $\dim D<\dim Z=3$, i.e. $\dim D=2$. Moreover, note that
$[D]_D\cong\mathcal{O}_D$. Thus by the adjunction formula we have
$(K_Z)_D\cong K_D-[D]_D\cong K_D$. Since $-K_D$ is ample, this
implies that $K_Z$ is not nef.

\smallskip

In cases $(12)$ and $(13)$, by similar arguments as in $(11)$, we
see that $K_Z$ cannot be nef.

\smallskip

Consider now case $(14)$. Let $\pi$ be the projection map $X\to
S$. Since $s_F\in\Gamma(\mathcal{E}|_F)$ for a general fibre
$F\cong\mathbb{P}^{n-2}$ of $\pi$, putting $D:=Z\cap F$ we have
$\dim D\geq n-2-r=\dim Z-2=1$. Since $D\subset F, D\subset Z$, we
deduce that
$$(K_Z)_D\cong (K_X+\det\mathcal{E})|_D\cong
(K_F+\det\mathcal{E}|_F)|_D = \mathcal{O}_F(-1)|_D.$$ This shows
that $K_Z$ cannot be nef.

\smallskip

In cases $(15)$ and $(16)$, by arguing as in case $(14)$, we
deduce that $K_Z$ is not nef.

\smallskip

Assume now that $(X,\mathcal{E})$ is as in case $(17)$. In this
situation, from section 2 we know that
$L:=(n-4)K_X+(n-3)\det\mathcal{E}$ is an ample line bundle on $X$
such that $L_F=\mathcal{O}_{\mathbb{P}^{n-3}}(1)$ for any fibre
$F\cong\mathbb{P}^{n-3}$ of $\pi:X\to V$. Since
$X\cong\mathbb{P}_V(\mathcal{F})$, we see that there exists an
ample vector bundle $\mathcal{V}$ of rank $n-2$ on $V$ such that
$(X,L)\cong (\mathbb{P}(\mathcal{V}),\xi_{\mathcal{V}})$, where
$\xi_{\mathcal{V}}$ is the tautological line bundle on
$\mathbb{P}(\mathcal{V})$. Thus
\begin{equation}\label{equation}
K_X+\det\mathcal{E}\cong
-(n-2)\xi_{\mathcal{V}}+\pi^*(K_V+\det\mathcal{V})+\det\mathcal{E}.
\end{equation}
By $\left[\det\mathcal{E}-(n-3)\xi_{\mathcal{V}}\right]|_F =
\mathcal{O}_F$, the relation (\ref{equation}) becomes
$$K_X+\det\mathcal{E}\cong
-\xi_{\mathcal{V}}+\pi^*(K_V+\det\mathcal{V}+\mathcal{H}),$$ for
some $\mathcal{H}\in\mathrm{ Pic }(V)$. Choose a very ample line
bundle $H$ on $V$ such that
$H':=K_V+\det\mathcal{V}+\mathcal{H}+H$ is ample on $V$. Put
$L':=L+\pi^*H$. Note that also $L'$ is ample and
\begin{equation}\label{equation-1}
K_X+\det\mathcal{E}+L'\cong \pi^*(H').
\end{equation}
Since $\dim Z = 3$, \cite[theorem 1.1]{LM1} tells us that
$\rho(X)=\rho(Z)$. Obviously, $\rho(X)=\rho(V)+1$. Hence
$\rho(Z)=\rho(V)+1$ and this shows that $Z$ cannot be isomorphic
to $V$ via $\pi$. Moreover, since
$s_F\in\Gamma(\mathcal{O}_{\mathbb{P}^{n-3}}(1)^{\oplus n-3})$ for
any fibre $F\cong\mathbb{P}^{n-3}$ of $\pi$, $D:=(s_F)_0=Z\cap F$
is a linear subspace of $F$, i.e. $D$ has degree one in
$F\cong\mathbb{P}^{n-3}$. Thus $\pi_Z:Z\to V$ is neither an
isomorphism nor a finite-to-one map onto $V$. So, there exists a
curve $\gamma\subset Z$ such that $\pi_Z(\gamma)$ is a point of
$V$. Restricting (\ref{equation-1}) to $Z$, we get
$K_Z+L_Z^{'}\cong\pi_{Z}^{*}H'$, and this gives
$(K_Z+L_Z^{'})\cdot\gamma = 0$, i.e. $K_Z\cdot\gamma =
-L_Z^{'}\cdot\gamma < 0$ by the ampleness of $L_Z^{'}$. This
implies that $K_Z$ is not nef.

\smallskip

Now, consider case $(18)$. Let $s_E$ denote the restriction of $s$
to $E$. Therefore
$s_E\in\Gamma(\mathcal{O}_{\mathbb{P}^{n-1}}(1)^{\oplus n-3})$.
Thus $D:=(s_E)_0=Z\cap E$ is a linear subspace of $E$ and we have
$$\dim D\geq n-1-(n-3)=2=\dim Z-1.$$ Obviously, $D\subseteq Z$. If
$D=Z$, then $Z\cong\mathbb{P}^3$, but this gives the contradiction
$\rho(X)=\rho(Z)=1$ by \cite[theorem 1.1]{LM1}. Therefore
$D\subsetneq Z$, and the irreducibility of $Z$ implies that $\dim
D < \dim Z$, i.e. $D\cong\mathbb{P}^2$. Furthermore,
$$\mathcal{O}_D(D)\cong\mathcal{O}_Z(Z\cap E)|_D\cong
\mathcal{O}_X(E)|_D\cong\mathcal{O}_E(E)|_D\cong\mathcal{O}_{\mathbb{P}^2}(-1).$$
Thus by the adjunction formula applied to $D\subset Z$, we get
$$(K_Z)_D\cong
K_D-\mathcal{O}_D(D)=\mathcal{O}_{\mathbb{P}^2}(-3)+\mathcal{O}_{\mathbb{P}^2}(1)=\mathcal{O}_{\mathbb{P}^2}(-2).$$
This shows that $K_Z$ is not nef.

\smallskip

Suppose that $(X,\mathcal{E})$ is as in case $(19)$. Let $s_{E_i}$
denote the restriction of $s$ to $E_i$. Then
$s_{E_i}\in\Gamma(\mathcal{E}|_{E_i})$. So $D:=(s_{E_i})_0=Z\cap
E_i$ is a linear subspace of $E_i$ and
$$\dim D\geq n-1-(n-3)=2=\dim Z-1.$$ By \cite[theorem 1.1]{LM1}
note that $D\neq Z$ since $\dim Z=3$. Thus $D\subsetneq Z$ and the
irreducibility of $Z$ implies that $D$ is $\mathbb{P}^2$.
Furthermore,
\begin{equation}
\mathcal{O}_D(D)\cong\mathcal{O}_Z(Z\cap
E_i)|_D\cong\mathcal{O}_{E_i}(E_i)|_D\cong
\mathcal{O}_{\mathbb{P}^2}(-2) \notag
\end{equation}
Thus by the adjunction formula applied to $D\subset Z$, we get
$(K_Z)|_D\cong K_D-\mathcal{O}_D(D)=\mathcal{O}_D(-1)$, i.e. $K_Z$
cannot be nef.

\smallskip

Finally, assume that we are in case $(20)$. Suppose that $D_i$ is
as in cases (j) and (jj). Let $s_{D_i}$ denote the restriction of
$s$ to $D_i$. Then $s_{D_i}\in\Gamma(\mathcal{E}|_{D_i})$. So
$D:=(s_{D_i})_0=Z\cap D_i$ is a quadric hypersurface contained in
$D_i$ and
$$\dim D\geq n-1-(n-3)=2=\dim Z-1.$$ By \cite[theorem 1.1]{LM1}
note that $D\neq Z$ since $\dim Z=3$. Thus $D\subsetneq Z$ and the
irreducibility of $Z$ implies that $D$ is a (possible singular)
quadric $\mathbb{Q}^2\subset\mathbb{P}^3$. Furthermore,
\begin{equation}
\mathcal{O}_D(D)\cong\mathcal{O}_Z(Z\cap
D_i)|_D\cong\mathcal{O}_{D_i}(D_i)|_D\cong
\mathcal{O}_{\mathbb{Q}^2}(-1) \notag
\end{equation}
Thus by the adjunction formula applied to $D\subset Z$, we have
$(K_Z)|_D\cong K_D-\mathcal{O}_D(D)=\mathcal{O}_D(-1)$. Therefore
$K_Z$ is not nef.

Suppose now that $D_i$ is as in case (jjj). Let $s_{F_i}$ be the
restriction of $s$ to any fiber $F_i\cong\mathbb{P}^{n-2}$ of
$D_i\to Z_i$. Then
$s_{F_i}\in\Gamma(\mathcal{O}_{\mathbb{P}^{n-2}}(1)^{\oplus
n-3})$. So $D:=(s_{F_i})_0=Z\cap F_i$ is a linear subspace of
$F_i$ and
$$\dim D\geq n-2-(n-3)=1=\dim Z-2.$$ Since $\rho(X)\neq 1$, by
\cite[theorem 1.1]{LM1} we see that $D\neq Z$. Thus the
irreducibility of $Z$ implies that $1\leq\dim D\leq 2$. Since
$D\subset Z$ and $D\subset F_i\subset D_i$, we obtain that
$$K_Z|_D\cong (K_X+\det\mathcal{E})|_D\cong
\left[(K_X+\det\mathcal{E})_{D_i}\right]|_D\cong
\left[K_{D_i}-\mathcal{O}_{D_i}(D_i)+\det\mathcal{E}_{D_i}\right]|_D\cong$$
$$\cong \left[K_{F}-\mathcal{O}_{D_i}(D_i)|_F+\det\mathcal{E}_{F}\right]|_D\cong \mathcal{O}_F(-1)|_D.$$
This shows that $K_Z$ cannot be nef. $\square$

\smallskip

When $\dim Z\geq 4$, having in mind \cite[theorems $1,2$ and
$3$]{LM}, we get the following

\begin{Prop}\label{dim Z>3}
Let $X$ be an $n$-fold. Let $\mathcal{E}$ be an ample vector
bundle of rank $r\geq 2$ on $X$ such that there exists a global
section $s\in\Gamma(\mathcal{E})$ whose zero locus $Z=(s)_0\subset
X$ is a smooth subvariety of dimension $n-r\geq 4$ and let $H$ be
an ample line bundle on $X$. Suppose that $K_Z+(\dim Z-2)H_Z$ is
nef. Then $K_Z+(\dim Z-3)H_Z$ is nef if and only if
$(X,\mathcal{E},H)$ is not any of the following triplets:
\begin{enumerate}
\item[$(1')$] $X\cong\mathbb{P}^n$ and $(\mathcal{E}, H)$ is one
of the following pairs:
\begin{enumerate}
\item[$(i)$] $(\mathcal{O}_{\mathbb{P}^n}(1)^{\oplus
n-4},\mathcal{O}_{\mathbb{P}^n}(4))$; \ $(ii)$\
$(\mathcal{O}_{\mathbb{P}^n}(1)^{\oplus
n-4},\mathcal{O}_{\mathbb{P}^n}(3))$; \item[$(iii)$]
$(\mathcal{O}_{\mathbb{P}^n}(1)^{\oplus
n-5}\oplus\mathcal{O}_{\mathbb{P}^n}(2),\mathcal{O}_{\mathbb{P}^n}(3))$;
\ $(iv)$\ $(\mathcal{O}_{\mathbb{P}^n}(1)^{\oplus
n-5}\oplus\mathcal{O}_{\mathbb{P}^n}(2),\mathcal{O}_{\mathbb{P}^n}(2))$;
\item[$(v)$] $(\mathcal{O}_{\mathbb{P}^n}(1)^{\oplus
n-6}\oplus\mathcal{O}_{\mathbb{P}^n}(2)^{\oplus
2},\mathcal{O}_{\mathbb{P}^n}(2))$; \ $(vi)$\
$(\mathcal{O}_{\mathbb{P}^n}(1)^{\oplus
n-5}\oplus\mathcal{O}_{\mathbb{P}^n}(3),\mathcal{O}_{\mathbb{P}^n}(2))$;
\item[$(vii)$] $(\mathcal{O}_{\mathbb{P}^n}(1)^{\oplus
n-5},\mathcal{O}_{\mathbb{P}^n}(2))$; \ $(viii)$\
$(\mathcal{O}_{\mathbb{P}^n}(1)^{\oplus
n-6}\oplus\mathcal{O}_{\mathbb{P}^n}(2),\mathcal{O}_{\mathbb{P}^n}(2))$;
\item[$(ix)$] $(\mathcal{O}_{\mathbb{P}^n}(1)^{\oplus
n-6},\mathcal{O}_{\mathbb{P}^n}(2))$; \ $(x)$\
$c_1(\mathcal{E})=r+3$ and $H=\mathcal{O}_{\mathbb{P}^n}(1)$;
\end{enumerate}
\item[$(2')$] $X\cong\mathbb{Q}^n$ and $(\mathcal{E}, H)$ is one
of the following pairs:
\begin{enumerate}
\item[(l)] $(\mathcal{O}_{\mathbb{Q}^n}(1)^{\oplus
n-4},\mathcal{O}_{\mathbb{Q}^n}(3))$; \item[(ll)]
$(\mathcal{O}_{\mathbb{Q}^n}(1)^{\oplus
n-4},\mathcal{O}_{\mathbb{Q}^n}(2))$; \item[(lll)]
$(\mathcal{O}_{\mathbb{Q}^n}(1)^{\oplus
n-5}\oplus\mathcal{O}_{\mathbb{Q}^n}(2),\mathcal{O}_{\mathbb{Q}^n}(2))$;
\item[(lv)] $(\mathcal{O}_{\mathbb{Q}^n}(1)^{\oplus
n-5},\mathcal{O}_{\mathbb{Q}^n}(2))$; \item[(v)]
$c_1(\mathcal{E})=\mathcal{O}_{\mathbb{Q}^n}(r+2)$ and
$H=\mathcal{O}_{\mathbb{Q}^n}(1)$;
\end{enumerate}
\item[$(3')$] $X$ is a Del Pezzo $n$-fold with $\mathrm{Pic
}(X)\cong\mathbb{Z}[\mathcal{O}_X(1)]$ and one of the following
possibilities can occur:
\begin{enumerate}
\item[$(\alpha)$] $\mathcal{E}\cong\mathcal{O}_X(1)^{\oplus r}$
and $H=\mathcal{O}_X(1)$; \item[$(\beta)$] for any line $l$ of
$(X,\mathcal{O}_X(1))$ either
$\mathcal{E}|_l\cong\mathcal{O}_l(1)^{\oplus n-4}$ and
$H|_l=\mathcal{O}_l(2)$, or
$\mathcal{E}|_l\cong\mathcal{O}_l(1)^{\oplus
r-1}\oplus\mathcal{O}_l(2)$ and $H|_l=\mathcal{O}_l(1)$;
\end{enumerate}
\item[$(4')$] $X$ is a Mukai $n$-fold with $\mathrm{Pic
}(X)\cong\mathbb{Z}[\mathcal{O}_X(1)]$ and
$\mathcal{E}\cong\mathcal{O}_X(1)^{\oplus r}$; \item[$(5')$]
$(\mathbb{P}^3\times\mathbb{P}^3,
\mathcal{O}_{\mathbb{P}^3\times\mathbb{P}^3}(1,1)^{\oplus 2},
\mathcal{O}_{\mathbb{P}^3\times\mathbb{P}^3}(1,1))$; \item[$(6')$]
there exists a vector bundle $\mathcal{V}$ on a smooth curve $C$
such that $X\cong\mathbb{P}_C(\mathcal{V})$; moreover, for any
fibre $F\cong\mathbb{P}^{n-1}$ of $X\to C$, we have
$(\mathcal{E}|_F, H_F)$ is isomorphic to one of the following
pairs:
\begin{enumerate}
\item[$(a)$] $(\mathcal{O}_{\mathbb{P}^{n-1}}(1)^{\oplus n-4},
\mathcal{O}_{\mathbb{P}^{n-1}}(2))$; \item[$(b)$]
$(\mathcal{O}_{\mathbb{P}^{n-1}}(2)\oplus\mathcal{O}_{\mathbb{P}^{n-1}}(1)^{\oplus
n-5}, \mathcal{O}_{\mathbb{P}^{n-1}}(2))$; \item[$(c)$]
$(\mathcal{O}_{\mathbb{P}^{n-1}}(1)^{\oplus n-5},
\mathcal{O}_{\mathbb{P}^{n-1}}(2))$; \item[$(d)$]
$(\mathcal{O}_{\mathbb{P}^{n-1}}(1)^{\oplus n-4},
\mathcal{O}_{\mathbb{P}^{n-1}}(3))$; \item[$(e)$]
$(\mathcal{O}_{\mathbb{P}^{n-1}}(3)\oplus\mathcal{O}_{\mathbb{P}^{n-1}}(1)^{\oplus
r-1}, \mathcal{O}_{\mathbb{P}^{n-1}}(1))$; \item[$(f)$]
$(\mathcal{O}_{\mathbb{P}^{n-1}}(2)^{\oplus
2}\oplus\mathcal{O}_{\mathbb{P}^{n-1}}(1)^{\oplus r-2},
\mathcal{O}_{\mathbb{P}^{n-1}}(1))$;
\end{enumerate}
\item[$(7')$] $X$ is a section of a divisor of relative degree two
in a projective space $\mathbb{P}_{C}(\mathcal{G})$, where
$\mathcal{G}$ is a vector bundle of rank $n+1$ on a smooth curve
$C$; moreover, for any smooth fibre $F\cong\mathbb{Q}^{n-1}$ of
$X\to C$, where $\mathbb{Q}^{n-1}$ is a smooth quadric
hypersurface of $\mathbb{P}^n$, we have $(\mathcal{E}|_F,H_F)$ is
isomorphic to either $(\mathcal{O}_{\mathbb{Q}^{n-1}}(1)^{\oplus
n-4}, \mathcal{O}_{\mathbb{Q}^{n-1}}(2))$, or
$(\mathcal{O}_{\mathbb{Q}^{n-1}}(2)\oplus\mathcal{O}_{\mathbb{Q}^{n-1}}(1)^{\oplus
r-1}, \mathcal{O}_{\mathbb{Q}^{n-1}}(1))$; \item[$(8')$] the map
$\Phi:X\to C$ associated to the linear system
$|(n-3)K_X+(n-2)\det\mathcal{E}|$ makes $X$ a Del Pezzo fibration
over a smooth curve $ C$; moreover, any general smooth fibre $F$
of $\Phi$ is a Del Pezzo $(n-1)$-fold with $\mathrm{Pic
}(F)\cong\mathbb{Z}[\mathcal{O}_F(1)]$ such that
$\mathcal{E}|_F\cong\mathcal{O}_F(1)^{\oplus r}$ and
$H_F\simeq\mathcal{O}_F(1)$; \item[$(9')$] there exists a vector
bundle $\mathcal{V}$ on a smooth surface $S$ such that
$X\cong\mathbb{P}_S(\mathcal{V})$; moreover, for any fibre
$F\cong\mathbb{P}^{n-2}$ of $X\to S$, we have
$(\mathcal{E}|_F,H_F)$ is isomorphic to either
$(\mathcal{O}_{\mathbb{P}^{n-2}}(1)^{\oplus n-4},
\mathcal{O}_{\mathbb{P}^{n-2}}(2))$, or
$(\mathcal{O}_{\mathbb{P}^{n-2}}(2)\oplus\mathcal{O}_{\mathbb{P}^{n-2}}(1)^{\oplus
r-1},\mathcal{O}_{\mathbb{P}^{n-2}}(1))$; \item[$(10')$] the map
$\psi:X\to S$ associated to the linear system
$|(n-3)K_X+(n-2)\det\mathcal{E}|$ makes $X$ a quadric fibration
over a smooth surface $S$; moreover, for any general fibre
$F\cong\mathbb{Q}^{n-2}$ we have $(\mathcal{E}|_F,H_F)\cong
(\mathcal{O}_{\mathbb{Q}^{n-2}}(1)^{\oplus
r},\mathcal{O}_{\mathbb{Q}^{n-2}}(1))$, where $\mathbb{Q}^{n-2}$
is a smooth quadric hypersurface of $\mathbb{P}^{n-1}$;
\item[$(11')$] there is a vector bundle $\mathcal{F}$ on a smooth
$3$-fold $V$ such that $X\cong\mathbb{P}_V(\mathcal{F})$;
moreover, for any fibre $F\cong\mathbb{P}^{n-3}$ of $X\to V$, we
get $(\mathcal{E}|_F,H_F)\cong
(\mathcal{O}_{\mathbb{P}^{n-3}}(1)^{\oplus r},
\mathcal{O}_{\mathbb{P}^{n-3}}(1))$; \item[$(12')$] the map $\psi
:X\to X'$ associated to the linear system
$|(n-2)K_X+(n-1)\det\mathcal{E}|$ is a birational morphism which
contracts an extremal face spanned by extremal rays $R_i$ for some
$i$ in a finite set of index. Let $\psi_i:X\to X_i$ be the
contraction associated to $R_i$. Then each $\psi_i$ is birational
and of divisorial type; moreover, if $E_i$ is an exceptional
divisor of $\psi_i$, then $E_i\cong\mathbb{P}^{n-1}$ and
$([E_i]_{E_i},\mathcal{E}|_{E_i}, H_{E_i})\cong
(\mathcal{O}_{\mathbb{P}^{n-1}}(-2),\mathcal{O}_{\mathbb{P}^{n-1}}(1)^{\oplus
r},\mathcal{O}_{\mathbb{P}^{n-1}}(1))$; \item[$(13')$] the map
$\phi:X\to X'$ associated to the linear system
$|(n-3)K_X+(n-2)\det\mathcal{E}|$ is a birational morphism which
contracts an extremal face. Let $R_i$ be the extremal rays
spanning this face for some $i$ in a finite set of index. Call
$\rho_i:X\to W_i$ the contraction associated to one of the $R_i$.
Then we have each $\rho_i$ is birational and of divisorial type;
if $D_i$ is one of the exceptional divisors and $Z_i=\rho_i(D_i)$,
we have $\mathrm{dim } Z_i\leq 1$ and one of the following
possibilities can occur:
\begin{enumerate}
\item[(j)] $\mathrm{dim } Z_i=0, D_i\cong\mathbb{P}^{n-1}$ and
$([D_i]_{D_i},\mathcal{E}|_{D_i},H_{D_i})$ is either
\begin{enumerate}
\item[$(j_1)$] $(\mathcal{O}_{\mathbb{P}^{n-1}}(-1),
\mathcal{O}_{\mathbb{P}^{n-1}}(1)^{\oplus n-4},
\mathcal{O}_{\mathbb{P}^{n-1}}(2))$; or \item[$(j_2)$]
$(\mathcal{O}_{\mathbb{P}^{n-1}}(-1),
\mathcal{O}_{\mathbb{P}^{n-1}}(2)\oplus\mathcal{O}_{\mathbb{P}^{n-1}}(1)^{\oplus
r-1}, \mathcal{O}_{\mathbb{P}^{n-1}}(1))$;
\end{enumerate}
\item[(jj)] $\mathrm{dim } Z_i=0, D_i$ is a (possible singular)
quadric hypersurface $\mathbb{Q}^{n-1}\subset\mathbb{P}^{n}$ and
$([D_i]_{D_i},\mathcal{E}|_{D_i},H_{D_i})\cong
(\mathcal{O}_{\mathbb{Q}^{n-1}}(-1),\mathcal{O}_{\mathbb{Q}^{n-1}}(1)^{\oplus
r},\mathcal{O}_{\mathbb{Q}^{n-1}}(1))$; \item[(jjj)] $\mathrm{dim
} Z_i=1, W_i$ and $Z_i$ are smooth projective varieties and
$\rho_i$ is the blow-up of $W_i$ along $Z_i$; furthermore, the
triplets $({\mathcal{O}_{D_i}(D_i)}|_{F_i}, \mathcal{E}|_{F_i},
H_{F_i})$ are isomorphic to $(\mathcal{O}_{\mathbb{P}^{n-2}}(-1),
\mathcal{O}_{\mathbb{P}^{n-2}}(1)^{\oplus r},
\mathcal{O}_{\mathbb{P}^{n-2}}(1))$ for any fibre
$F_i\cong\mathbb{P}^{n-2}$ of the restriction map $\rho_i|_{D_i}:
D_i\to Z_i$.
\end{enumerate}
Moreover, the map $\phi$ is a composition of disjoint extremal
contractions as in $\mathrm{(j), (jj)}$ and $\mathrm{(jjj)}$.
\end{enumerate}
\end{Prop}

\noindent\textbf{Proof.} To prove the `if' part, suppose that
$K_Z+(\dim Z-3)H_Z$ is not nef. Set
$\mathcal{E}':=\mathcal{E}\oplus H^{\oplus n-r-3}$. Then
$\mathcal{E}'$ is an ample vector bundle of rank $n-3$ on $X$.
Now, we get
$$(K_{X}+\det\mathcal{E}')|_Z=[K_X+\det\mathcal{E}+(\dim
Z-3)H]|_Z\cong K_Z+(\dim Z-3)H_Z .$$ So, $K_X+\det\mathcal{E}'$ is
not nef on $X$. Thus we deduce that $(X,\mathcal{E}')$ is as in
Theorem \ref{thm}. Keeping in mind that
$\mathcal{E}'=\mathcal{E}\oplus H^{\oplus n-r-3}$ with $n\geq
r+4\geq 6$, by \cite[theorems $1,2$ and $3$]{LM} we can easily see
that cases $(8)$, $(9)$, $(10)$, $(11)(a)$, $(12)(i)(iii)$, $(18)$
of Theorem \ref{thm} cannot occur. Moreover, as to case $(15)$, we
have $H|_F\cong\mathcal{O}_{\mathbb{P}^{n-2}}(1)$ for any closed
fiber $F\cong\mathbb{P}^{n-2}$ of $\pi:X\to S$. This shows that
$(X,H)$ is an adjunction-theoretic scroll over $S$ and by
\cite[proposition 3.2.1 and theorem 14.1.1]{BS} and
\cite[(2.12)]{Fujita} we can conclude that there exists a suitable
vector bundle $\mathcal{V}$ of rank $n-1$ on $S$ such that
$(X,H)\cong (\mathbb{P}_S(\mathcal{V}),\xi_{\mathcal{V}})$, where
$\xi_{\mathcal{V}}$ is the tautological line bundle on
$\mathbb{P}_S(\mathcal{V})$. From \cite[theorem 3(14)]{LM} it
follows that $K_Z+(\dim Z-2)H_Z$ cannot be nef, a contradiction.
In all the other cases, by \cite[theorems $1,2$ and $3$]{LM} and
easy computations we are done.

\medskip

\noindent Conversely, if $(X,\mathcal{E}, H)$ is as in cases
$(1')$ to $(5')$, by the adjunction formula we see that $K_Z+(\dim
Z-3)H_Z$ is not nef.

\smallskip

Suppose we are in case $(6')$. Let $\pi$ be the projection map
$X\to C$ and let $s_F$ denote the restriction of the global
section $s$ to a general fibre $F\cong\mathbb{P}^{n-1}$ of $\pi$.
Then $s_F\in\Gamma (\mathcal{E}|F)$ and putting $D:=(s_F)_0=Z\cap
F$, we have $\dim D\geq n-r-1\geq 3$. On the other hand, since $Z$
is irreducible, we have $\dim D<\dim Z=n-r$, i.e. $\dim D=\dim
Z-1=n-r-1$. Moreover, note that $[D]_D\cong\mathcal{O}_D$. Thus
$(K_Z)_D\cong K_D-[D]_D\cong K_D$. Hence $-[K_Z+(\dim
Z-3)H_Z]|_D=-[K_D+(\dim D-2)H_D]$ is ample on $D\subset Z$ and
this shows that $K_Z+(\dim Z-3)H_Z$ cannot be nef.

\smallskip

In cases $(7')$ and $(8')$, by a similar argument as in $(5')$, we
can see that $K_Z+(\dim Z-3)H_Z$ is not nef.

\smallskip

Consider now case $(9')$. Let $\pi$ be the projection map $X\to
S$. Since $s_F\in\Gamma (\mathcal{E}|F)$ for a general fibre
$F\cong\mathbb{P}^{n-2}$ of $\pi$, putting $D:=Z\cap F=(s_F)_0$ we
get $\dim D\geq n-r-2=\dim Z-2\geq 2$. Since $D\subset F$ and
$D\subset Z$, we deduce in both situations that
$$(K_Z+(\dim Z-3)H_Z)|_D=(K_X+\det\mathcal{E}+(\dim Z-3)H)|_D=$$
$$=[K_F+\det\mathcal{E}|_F+(\dim
Z-3)H_F]|_D\cong\mathcal{O}_{\mathbb{P}^{n-2}}(-1)|_D.$$ This
implies that $K_Z+(\dim Z-3)H_Z$ cannot be nef.

\smallskip

By arguing as in $(9')$, we can see that also in case $(10')$ the
line bundle $K_Z+(\dim Z-3)H_Z$ is not nef.

\smallskip

Assume now that $(X,\mathcal{E},H)$ is as in $(11')$. Let $\pi$ be
the projection map $\mathbb{P}_V(\mathcal{V})\to V$. Since
$s_F\in\Gamma (\mathcal{E}|_F)$ for any fibre
$F\cong\mathbb{P}^{n-3}$ of $\pi$, $D:=(s_F)_0=Z\cap F$ is not
empty for any fibre $F$, i.e. the restriction map $\pi_Z:Z\to V$
of $\pi$ to $Z$ is surjective. On the other hand, we have $\dim
Z>\dim D\geq n-r-3=\dim Z-3$, because of the irreducibility of
$Z$. So $\dim Z-3\leq\dim D\leq\dim Z-1$. Furthermore, since
$K_X\cong -(n-2)\xi_{\mathcal{V}}+\pi^*(K_V+\det\mathcal{V})$,
where $\xi_{\mathcal{V}}$ is the tautological line bundle on
$\mathbb{P}_V(\mathcal{V})$, we get
\begin{equation}\label{*}
K_X+\det\mathcal{E}+(n-r-3)H=[\det\mathcal{E}+(n-r-3)H-(n-3)\xi_{\mathcal{V}}]
+\pi^*(K_V+\det\mathcal{V})-\xi_{\mathcal{V}}.
\end{equation}
Since
$[\det\mathcal{E}+(n-r-3)H-(n-3)\xi_{\mathcal{V}}]|_F\cong\mathcal{O}_F$,
it follows that
$$\det\mathcal{E}+(n-r-3)H-(n-3)\xi_{\mathcal{V}}\cong\pi^*\mathcal{D},$$
for some $\mathcal{D}\in\mathrm{Pic } (V)$. By (\ref{*}) we obtain
that
$$[K_Z+(\dim Z-3)H_Z]|_D=([K_X+\det\mathcal{E}+(n-r-3)H]_Z)|_D\cong $$
$$\cong
([K_X+\det\mathcal{E}+(n-r-3)H]_F)|_D=-({\xi_{\mathcal{V}}}_F)|_D\cong\mathcal{O}_{\mathbb{P}^{n-3}}(-1)|_D,$$
for a general fibre $D=F|_Z$ of $\pi_Z$. This implies that
$K_Z+(\dim Z-3)H_Z$ is not nef.

\smallskip

Suppose we are in case $(12')$. Let $s_{E_i}$ denote the
restriction of the section $s$ to $E_i$. Then $s_{E_i}\in\Gamma
(\mathcal{E}|_{E_i})$ and so $D:=(s_{E_i})_0=E_i\cap Z$ is a
linear subspace of $E_i\cong\mathbb{P}^{n-1}$. Moreover, $\dim
D\geq n-r-1=\dim Z-1\geq 3$. Obviously, $D\subseteq Z$. If $D=Z$,
then $Z\subset E_i\subset X$ and we have the canonical surjection
$\mathcal{N}_{Z/X}\to\mathcal{N}_{E_i/X}|_Z$ between normal
bundles. Since $\mathcal{N}_{Z/X}\cong\mathcal{E}|_Z$ is ample, we
deduce that $\mathcal{N}_{E_i/X}|_Z$ is ample, but this gives a
contradiction since $E_i$ is an exceptional divisor on $X$. Thus
$D\subsetneqq Z$ and the irreducibility of $Z$ implies that $\dim
D<\dim Z$, i.e. $\dim D=\dim Z-1=n-r-1$. Therefore $D$ is
isomorphic to $\mathbb{P}^{n-r-1}$ and

\begin{equation}\label{**''}
H_Z=H_D=(H_{E_i})|_D=\mathcal{O}_{\mathbb{P}^{n-r}}(1) \notag
\end{equation}

\noindent Note that
\begin{equation}
\mathcal{O}_D(D)\cong\mathcal{O}_Z(Z\cap
E_i)|_D\cong\mathcal{O}_X(E_i)|_D\cong\mathcal{O}_{E_i}(E_i)|_D\cong\mathcal{O}_{\mathbb{P}^{n-r-1}}(-2)
\notag
\end{equation}
and
\begin{equation}
(D,H_D,\mathcal{O}_D(D))\cong (\mathbb{P}^{n-r-1},
\mathcal{O}_{\mathbb{P}^{n-r-1}}(1),
\mathcal{O}_{\mathbb{P}^{n-r-1}}(-2)). \notag
\end{equation}
Thus by the adjunction formula applied to $D\subset Z$, we obtain
that
\begin{equation}
(K_Z)|_D\cong
K_D-\mathcal{O}_D(D)\cong\mathcal{O}_{\mathbb{P}^{n-r-1}}(-n+r+2).
\notag
\end{equation}
Hence $(K_Z+(\dim Z-3)H_Z)|_D\cong\mathcal{O}_{D}(-1)$ and this
shows that $K_Z+(\dim Z-3)H_Z$ cannot be nef.

\smallskip

Consider now case $(13')$(j). Let $s_{D_i}$ denote the restriction
of the section $s$ to $D_i$. Then $s_{D_i}\in\Gamma
(\mathcal{E}|_{D_i})$ and so $D:=(s_{D_i})_0=D_i\cap Z$ is either
a linear subspace of $D_i\cong\mathbb{P}^{n-1}$, or a (possible
singular) quadric hypersurface on $D_i$. Moreover, $\dim D\geq
n-r-1=\dim Z-1\geq 3$. Obviously, $D\subseteq Z$. If $D=Z$, then
$Z\subset D_i\subset X$ and we have the canonical surjection
$\mathcal{N}_{Z/X}\to\mathcal{N}_{D_i/X}|_Z$ between normal
bundles. Since $\mathcal{N}_{Z/X}\cong\mathcal{E}|_Z$ is ample, we
deduce that $\mathcal{N}_{D_i/X}|_Z$ is ample, but this gives a
contradiction since $D_i$ is an exceptional divisor on $X$.
Thus $D\subsetneqq Z$ and the irreducibility of $Z$ implies that
$\dim D<\dim Z$, i.e. $\dim D=\dim Z-1=n-r-1$. Furthermore, $D$ is
isomorphic to either $\mathbb{P}^{n-r-1}$, or
$\mathbb{Q}^{n-r-1}\subset\mathbb{P}^{n-r}$ and

\begin{equation}\label{**}
H_Z=H_D=(H_{D_i})|_D=\left\{
\begin{array}{ll}
\mathcal{O}_{\mathbb{P}^{4}}(2) & \ \ \ \mathrm{in\ case\ } (j_1) \\
\mathcal{O}_{\mathbb{Q}^{n-r}}(1) & \ \ \ \mathrm{in\ case\ }
(j_2)
\end{array}
\right.
\end{equation}

\noindent Note that
\begin{equation}
\mathcal{O}_D(D)\cong\mathcal{O}_Z(Z\cap
D_i)|_D\cong\mathcal{O}_X(D_i)|_D\cong\mathcal{O}_{D_i}(D_i)|_D\cong\left\{
\begin{array}{ll}
\mathcal{O}_{\mathbb{P}^{3}}(-1) & \ \ \ \mathrm{in\ case\ } (j_1) \\
\mathcal{O}_{\mathbb{Q}^{n-r-1}}(-1) & \ \ \ \mathrm{in\ case\ }
(j_2)
\end{array}
\right. . \notag
\end{equation}
Thus
\begin{equation}
(D,H_D,\mathcal{O}_D(D))\cong\left\{
\begin{array}{ll}
(\mathbb{P}^{3}, \mathcal{O}_{\mathbb{P}^{3}}(2), \mathcal{O}_{\mathbb{P}^{3}}(-1)) & \ \ \ \mathrm{in\ case\ } (j_1) \\
(\mathbb{Q}^{n-r-1}, \mathcal{O}_{\mathbb{Q}^{n-r-1}}(1),
\mathcal{O}_{\mathbb{Q}^{n-r-1}}(-1)) & \ \ \ \mathrm{in\ case\ }
(j_2)
\end{array}
\right. \notag
\end{equation}
and by the adjunction formula applied to $D\subset Z$, we obtain
that
\begin{equation}
(K_Z)|_D\cong K_D-\mathcal{O}_D(D)\cong\left\{
\begin{array}{ll}
\mathcal{O}_{\mathbb{P}^{3}}(-3) & \ \ \ \mathrm{in\ case\ } (j_1) \\
\mathcal{O}_{\mathbb{Q}^{n-r-1}}(-n+r+2) & \ \ \ \mathrm{in\ case\
} (j_2)
\end{array}
\right. . \notag
\end{equation}
Hence $(K_Z+(\dim Z-3)H_Z)|_D\cong\mathcal{O}_{D}(-1)$ and this
shows that $K_Z+(\dim Z-3)H_Z$ is not nef.

\smallskip

In case $(13')$(jj), by arguing as in $(11')$(j), we see that
$D:=Z\cap D_i\subsetneqq Z$ and $\dim D=n-r-1$. Hence
$D\cong\mathbb{Q}^{n-r-1}\subset D_i$ and $H_Z|_D\cong
(H_{D_i})|_D\cong\mathcal{O}_D(1)$. Moreover,
$\mathcal{O}_D(D)\cong\mathcal{O}_Z(Z\cap
D_i)|_D\cong\mathcal{O}_{D_i}(D_i)|_D\cong\mathcal{O}_D(-1)$ and
so $K_Z|_D\cong
K_D-\mathcal{O}_D(D)\cong\mathcal{O}_{\mathbb{Q}^{n-r-1}}(-n+r+2)$.
This gives $(K_Z+(\dim Z-3)H_Z)|_D\cong\mathcal{O}_D(-1)$, i.e.
$K_Z+(\dim Z-3)H_Z$ cannot be nef.

Finally, let us consider case $(13')$(jjj). Let $s_F$ denote the
restriction of the section $s$ to a fibre $F\cong\mathbb{P}^{n-2}$
of $D_i\to Z_i$. Then $s_F\in\Gamma
(\mathcal{O}_{\mathbb{P}^{n-2}}(1)^{\oplus r})$ with $r\leq n-4$.
Thus $D:=(s_F)_0=Z\cap F$ is a linear subspace of
$F\cong\mathbb{P}^{n-2}$ and $\dim D\geq n-r-2=\dim Z-2\geq 2$.
Obviously, $D\subseteq Z$. If $D=Z$, then $Z\cong\mathbb{P}^{n-r}$
and $H_Z=H_D=(H_F)|_D\cong\mathcal{O}_{\mathbb{P}^{n-r}}(1)$. Thus
$K_Z+(\dim Z-2)H_Z\cong\mathcal{O}_{\mathbb{P}^{n-r}}(-3)$, but
this is a contradiction. Therefore, $D\subsetneqq Z$ and the
irreducibility of $Z$ implies that $n-r-2\leq\dim D\leq n-r-1$.
Then
$$(K_Z+(\dim Z-3)H_Z)|_D\cong
[(K_X+\det\mathcal{E}+(n-r-3)H)|_F]|_D=$$
$$=[(K_X+\det\mathcal{E}+(n-r-3)H)|_{D_i}]|_D=[(K_{D_i}-[D_i]_{D_i}+\det\mathcal{E}|_{D_i}+(n-r-3)H_{D_i})|_F]|_D=$$
$$=[(K_F-\mathcal{O}_{D_i}(D_i)|_F+\det\mathcal{E}|_F+(n-r-3)H_F)]|_D=\mathcal{O}_{\mathbb{P}^{n-r}}(-1)|_D.$$
This shows that $K_Z+(\dim Z-3)H_Z$ cannot be nef again. $\square$

\medskip

\subsection{Some remarks on classical scrolls over fivefolds.}\label{3.2}

In very classical times scrolls naturally occurred very often. Let
us recall that by a \textit{classical scroll} we mean a
$\mathbb{P}^k$-bundle $X$ over a variety $Y$ together with an
ample line bundle $L$ such that
$L|_F\simeq\mathcal{O}_{\mathbb{P}^k}(1)$ for any fibre
$F\cong\mathbb{P}^k$ with $k=\mathrm{dim } X-\mathrm{dim } Y$.
This is equivalent to saying that $(X,L)\cong
(\mathbb{P}(\mathcal{E}),\xi_{\mathcal{E}})$, where
$\mathcal{E}=p_*L$ is an ample vector bundle of rank $k+1$ on $Y$
and $p:X\to Y$ is the projection map. In this situation, the
canonical bundle formula gives
\begin{equation}
K_X+(\mathrm{dim } X-\mathrm{dim } Y+1)L=K_X+(\mathrm{rk}\
\mathcal{E}) L\simeq p^*(K_Y+\det\mathcal{E}). \label{a4}
\end{equation}
From the adjunction theoretic point of view the correct definition
of scroll (see \cite{S}) is that of \textit{adjunction-theoretic
scroll} over a normal variety $Y$. This means that there exists a
morphism with connected fibres, $p:X\to Y$, such that
$$K_X+(\mathrm{dim } X-\mathrm{dim } Y+1)L\simeq p^* H,$$ for some ample line bundle
$H$ on $Y$. The general fibre is
$(\mathbb{P}^k,\mathcal{O}_{\mathbb{P}^k}(1))$ with $k=\mathrm{dim
} X-\mathrm{dim } Y$, but the special fibres can vary quite a lot.

Note that if $L$ is further very ample and $(X,L)$ is an
adjunction theoretic scroll over a normal variety $Y$ of dimension
$m\leq 4$, then $Y$ is smooth and $(X,L)$ is a
$\mathbb{P}^k$-bundle over $Y$ with $k=\mathrm{dim } X-m$. This
follows from a general result due to Sommese (\cite[theorem
3.3]{S}) for $m\leq 2$ (and $L$ merely spanned), from
\cite[proposition 3.2.3]{BSW} and \cite[proposition 2.1]{T} for
$m=3$ (and $L$ ample and spanned), and from \cite[theorem 2.2]{T}
for $m=4$.

The following results are concerned with the other direction. In
the stable range $\mathrm{dim } X\geq 2\mathrm{dim } Y-1$, a
classical scroll is also an adjunction theoretic scroll with a few
easy exceptions. These results depend on Mori theory \cite{Mori}
and \cite{F1}. In the unstable range $\mathrm{dim } X\leq
2\mathrm{dim } Y-2$ and $\mathrm{dim } Y\leq 4$, classical scrolls
that are not adjunction theoretic scrolls in the modern sense have
been classified in \cite{BS2} and \cite[$\S 3$]{T}.

As to the cases $n=2m-3\geq 7$, i.e. $k=m-3\geq 2$ and
$\mathrm{rk}\ \mathcal{E}=m-2$, let us give here this immediate
consequence of \cite[theorem 5.1]{AM}, \cite{Maeda} and
\cite[theorem 1]{Occhetta}, as noted in \cite[remark 3.3]{T}.

\begin{Prop}\label{class-scroll k=m-3}
Let $X$ be an $n$-fold and let $L$ be an ample line bundle on $X$.
Assume that $(X,L)$ $\cong$ $(\mathbb{P}
(\mathcal{E}),\mathcal{O}_{\mathbb{P} (\mathcal{E})}(1))$ is a
$\mathbb{P}^{n-m}$-bundle, $\pi : X\to Y$, over a smooth variety
$Y$ of dimension $m$ with $\mathcal{E}=\pi_*L$. If $n=2m-3\geq 7$,
then $(X,L)$ is an adjunction-theoretic scroll over $Y$ under
$\pi$ unless either:

\begin{enumerate}
    \item[$(1)$] $Y\cong\mathbb{P}^m$ and $\mathcal{E}\cong
    \mathcal{O}_{\mathbb{P}^m}(1)^{\oplus m-2}$,
    $\mathcal{O}_{\mathbb{P}^m}(2)\oplus\mathcal{O}_{\mathbb{P}^m}(1)^{\oplus m-3}$,
    $\mathcal{O}_{\mathbb{P}^m}(2)^{\oplus 2}\oplus\mathcal{O}_{\mathbb{P}^m}(1)^{\oplus m-4}$,
    $\mathcal{O}_{\mathbb{P}^m}(3)\oplus\mathcal{O}_{\mathbb{P}^m}(1)^{\oplus m-3}$;

    \item[$(2)$] $Y\cong\mathbb{Q}^m$ and $\mathcal{E}\cong\mathcal{O}_{\mathbb{Q}^m}(1)^{\oplus m-2}$,
    $\mathcal{O}_{\mathbb{Q}^m}(2)\oplus\mathcal{O}_{\mathbb{Q}^m}(1)^{\oplus m-3}$, where $\mathbb{Q}^m$ is a smooth quadric hypersurface of $\mathbb{P}^{m+1}$;

    \item[$(3)$] $Y$ is a Del Pezzo $m$-fold with $b_2(Y)=1$, i.e. $\mathrm{Pic}(Y)$ is generated by an ample line bundle
    $\mathcal{O}_Y(1)$ such that $-K_Y\cong (m-1)\mathcal{O}_Y(1)$
    and $\mathcal{E}\cong\mathcal{O}_Y(1)^{\oplus m-2}$;

    \item[$(4)$] there is a vector bundle $\mathcal{V}$ on a smooth curve $C$ such that $Y\cong\mathbb{P}(\mathcal{V})$
    and $\mathcal{E}|_F\cong\mathcal{O}_{\mathbb{P}^{m-1}}(1)^{\oplus m-2}$ for any fibre $F\cong\mathbb{P}^{m-1}$ of $Y\to C$;

    \item[$(5)$] there is a vector bundle $\mathcal{V}$ on a smooth curve $C$ such that $Y\cong\mathbb{P}(\mathcal{V})$
    and $\mathcal{E}|_F\cong\mathcal{O}_{\mathbb{P}^{m-1}}(2)\oplus\mathcal{O}_{\mathbb{P}^{m-1}}(1)^{\oplus m-3}$ for any fibre $F\cong\mathbb{P}^{m-1}$ of $Y\to C$;

    \item[$(6)$] there is a surjective morphism $q:Y\to C$ onto a smooth curve $C$ such that any general fibre $F$ of $q$ is a smooth
    quadric hypersurface $\mathbb{Q}^{m-1}\subset\mathbb{P}^{m}$ with $\mathcal{E}|_F\cong\mathcal{O}_{\mathbb{Q}^{m-1}}(1)^{\oplus m-2}$;

    \item[$(7)$] there is a vector bundle $\mathcal{V}$ on a smooth surface $S$ such that $Y\cong\mathbb{P}(\mathcal{V})$
    and $\mathcal{E}|_F\cong\mathcal{O}_{\mathbb{P}^{m-2}}(1)^{\oplus m-2}$ for any fibre $F\cong\mathbb{P}^{m-2}$ of $Y\to S$;

    \item[$(8)$] $Y$ is a Fano $m$-fold with
    $-K_Y\simeq\det\mathcal{E}$; moreover, if $m\geq 6$ then
    $b_2(Y)=1$ except for $Y\cong\mathbb{P}^3\times\mathbb{P}^3$
    and $\mathcal{E}\cong\mathcal{O}_{\mathbb{P}^3\times\mathbb{P}^3}(1,1)^{\oplus
    4}$;

    \item[$(9)$] there exists a morphism $\Phi:Y\to W$ onto a normal
    variety $W$ supported by (a large multiple of) $K_Y+\det\mathcal{E}$ and
    $\mathrm{dim } W\leq 3$; if $F$ is a general fibre of $\Phi$,
    then we have the following possibilities:

    \begin{enumerate}
    \item[(a)] $W$ is a smooth curve and $\Phi$ is flat (equidimensional)
    map; moreover, the pair $(F,\mathcal{E}|_F)$ is one of the
    following:
    \begin{enumerate}
    \item[(a1)] $F\cong\mathbb{P}^{m-1}$ and
    $\mathcal{E}|_F\cong\mathcal{O}(2)^{\oplus
    2}\oplus\mathcal{O}(1)^{m-4},
    \mathcal{O}(3)\oplus\mathcal{O}(1)^{m-3}$;
    \item[(a2)] $F\cong\mathbb{Q}^{m-1}$ and
    $\mathcal{E}|_F\cong\mathcal{O}(2)\oplus\mathcal{O}(1)^{\oplus
    m-3}$, where $\mathbb{Q}^{m-1}$ is the
    quadric hypersurface of $\mathbb{P}^{m}$;
    \item[(a3)] $F\cong\mathbb{Q}^{4}$ and $\mathcal{E}|_F\cong
    E(2)\oplus\mathcal{O}(1)$, where $E$ is a spinor bundle on the
    quadric hypersurface $\mathbb{Q}^4$ of $\mathbb{P}^5$; \item[(a4)] $F$
    is a Del Pezzo $(m-1)$-fold with $b_2(F)=1$, i.e.
    $\mathrm{Pic}(F)$ is generated by an ample line bundle
    $\mathcal{O}_F(1)$ such that $-K_F\simeq
    (m-2)\mathcal{O}_F(1)$ and
    $\mathcal{E}|_F\cong\mathcal{O}_F(1)^{\oplus m-2}$;
    \item[(a5)] $F\cong\mathbb{P}^2\times\mathbb{P}^2$ and
    $\mathcal{E}|_F\cong\mathcal{O}(1,1)^{\oplus 3}$.
    \end{enumerate}
    In particular, if $m\geq 6$ then $\Phi$ is an elementary contraction. Furthermore, if $F\cong\mathbb{P}^{m-1}$ then $Y$ is a classical
    scroll, while if $F\cong\mathbb{Q}^{m-1}\subset\mathbb{P}^{m}$ then $Y$ is a quadric bundle;
    \item[(b)] $W$ is a smooth surface, $\Phi$ is flat and for a general fibre $F$ of $\Phi$ the pair
    $(F,\mathcal{E}|_F)$ is one of the following:
    \begin{enumerate}
    \item[(b1)] $(\mathbb{P}^{m-2}, \mathcal{O}(2)\oplus\mathcal{O}(1)^{\oplus m-3})$;
    \item[(b2)] $(\mathbb{P}^{m-2}, T_{\mathbb{P}^{m-2}})$, where
    $T_{\mathbb{P}^{m-2}}$ is the tangent bundle on
    $\mathbb{P}^{m-2}$;
    \item[(b3)] $(\mathbb{Q}^{m-2}, \mathcal{O}(1)^{\oplus
    m-2})$, where $\mathbb{Q}^{m-2}$ is a smooth quadric hypersurface in
    $\mathbb{P}^{m-1}$.
    \end{enumerate}
    In particular, if $F\cong\mathbb{P}^{m-2}$ then all the fibres
    of $\Phi$ are $\mathbb{P}^{m-2}$;
    \item[(c)] $W$ is a $3$-fold and $\mathcal{E}|_F\cong\mathcal{O}(1)^{\oplus
    m-2}$ for all the fibres $F\cong\mathbb{P}^{m-3}$ of $\Phi$;
    \end{enumerate}
    \item[$(10)$] there exist an $m$-fold $Y'$, a
    morphism $\psi:Y\to Y'$ expressing $Y$ as blown up at a finite
    set $B$ of points and an ample vector bundle $\mathcal{E}'$ on $Y'$ such
    that $\mathcal{E}=\pi^*\mathcal{E}'\otimes [-\psi^{-1}(B)]$;
    or

    \item[$(11)$] a high multiple of $K_Y+\det\mathcal{E}$ defines
    a birational map, $\varphi: Y\to\widehat{Y}$, which contracts
    an extremal face. Let $R_i$, for $i$ in a finite set of index,
    the extremal rays spanning this face; call
    $\rho_i:Y\to\widehat{Y}_i$ the contraction associated to one
    of the $R_i$. Then we have each $\rho_i$ is birational
    and divisorial; if $D_i$ is one of the exceptional divisors
    and $Z_i=\rho_i(D_i)$, we have $\mathrm{dim } Z_i\leq 1$ and the
    following possibilities can occur:
    \begin{enumerate}
    \item[(i)] $\mathrm{dim } Z_i=0, D_i\cong\mathbb{P}^{m-1}$ and
    $[D_i]_{D_i}\simeq\mathcal{O}(-2)$, or $\mathcal{O}(-1)$;
    moreover, $\mathcal{E}|_{D_i}\cong\mathcal{O}(1)^{\oplus
    m-2}$, or $\mathcal{E}|_{D_i}\cong\mathcal{O}(2)\oplus\mathcal{O}(1)^{\oplus
    m-3}$;
    \item[(ii)] $\mathrm{dim } Z_i=0, D_i$ is a (possible singular)
    quadric hypersurface $\mathbb{Q}^{m-1}$ and
    $[D_i]_{D_i}\simeq\mathcal{O}(-1)$ with
    $\mathcal{E}|_{D_i}\cong\mathcal{O}(1)^{m-2}$;
    \item[(iii)] $\mathrm{dim } Z_i=1$, $Z_i$ and $\widehat{Y}_i$ are
    smooth projective varieties and $\rho_i$ is the blow-up of
    $\widehat{Y}_i$ along $Z_i$; moreover,
    $\mathcal{E}|_{f}\cong\mathcal{O}(1)^{\oplus m-2}$ for any
    fibre $f\cong\mathbb{P}^{m-2}$ of ${\rho_i}|_{D_i}:D_i\to
    Z_i$.
    \end{enumerate}
\end{enumerate}
\end{Prop}

Specializing the above result to the case $m=5$, we get also the
following

\begin{Cor}
Let $X$ be a smooth $n$-fold with $n\geq 7$ and $L$ an ample line
bundle on $X$. Assume that $(X,L)$ $\cong$ $(\mathbb{P}
(\mathcal{E}),\mathcal{O}_{\mathbb{P} (\mathcal{E})}(1))$ is a
$\mathbb{P}^k$-bundle, $\pi:X\to Y$, over a smooth $5$-fold $Y$
with $\mathcal{E}=\pi_*L$. Then $(X,L)$ is an adjunction scroll
over $Y$ under $\pi$ unless either
\begin{enumerate}

\item[$(1')$] $n=10$, $Y\cong\mathbb{P}^5$ and
$\mathcal{E}\cong\mathcal{O}_{\mathbb{P}^5}(1)^{\oplus 6}$;

\item[$(2')$] $n=9$, $Y\cong\mathbb{P}^5$ and
$\mathcal{E}\cong\mathcal{O}_{\mathbb{P}^5}(1)^{\oplus 5},
\mathcal{O}_{\mathbb{P}^5}(2)\oplus\mathcal{O}_{\mathbb{P}^5}(1)^{\oplus
4}$ or $T_{\mathbb{P}^5}$, where $T_{\mathbb{P}^5}$ is the tangent
bundle of $\mathbb{P}^5$;

\item[$(3')$] $n=9$, $Y\cong\mathbb{Q}^5$ and
$\mathcal{E}\cong\mathcal{O}_{\mathbb{Q}^5}(1)^{\oplus 5}$, where
$\mathbb{Q}^5$ is a smooth quadric hypersurface of $\mathbb{P}^6$;

\item[$(4')$] $n=9$, there is a rank $5$ vector bundle
$\mathcal{V}$ over a smooth curve $C$ such that
$Y\cong\mathbb{P}(\mathcal{V})$ and
$\mathcal{E}|_{f}\cong\mathcal{O}_f(1)^{\oplus 5}$ for any fibre
$f\cong\mathbb{P}^{4}$ of $Y\to C$;

\item[$(5')$] $n=8$, $Y\cong\mathbb{P}^5$ and $\mathcal{E}\cong
\mathcal{O}_{\mathbb{P}^5}(1)^{\oplus 4}$,
$\mathcal{O}_{\mathbb{P}^5}(2)\oplus\mathcal{O}_{\mathbb{P}^5}(1)^{\oplus
3}$, $\mathcal{O}_{\mathbb{P}^5}(2)^{\oplus
2}\oplus\mathcal{O}_{\mathbb{P}^5}(1)^{\oplus 2}$,
$\mathcal{O}_{\mathbb{P}^5}(3)\oplus\mathcal{O}_{\mathbb{P}^5}(1)^{\oplus
3}$;

\item[$(6')$] $n=8$, $Y\cong\mathbb{Q}^5$ and $\mathcal{E}$ is
either $\mathcal{O}_{\mathbb{Q}^5}(1)^{\oplus 4}$ or
$\mathcal{O}_{\mathbb{Q}^5}(2)\oplus\mathcal{O}_{\mathbb{Q}^5}(1)^{\oplus
3}$, where $\mathbb{Q}^5$ is a smooth quadric hypersurface of
$\mathbb{P}^{6}$;

\item[$(7')$] $n=8$, $Y$ is a Del Pezzo $5$-fold with $b_2(Y)=1$,
i.e. $\mathrm{Pic }(Y)$ is generated by an ample line bundle
$\mathcal{O}_Y(1)$ such that $-K_Y\simeq\mathcal{O}_Y(4)$ and
$\mathcal{E}\cong\mathcal{O}_Y(1)^{\oplus 4}$;

\item[$(8')$] $n=8$, there is a vector bundle $\mathcal{V}$ on a
smooth curve $C$ such that $Y\cong\mathbb{P}_C(\mathcal{V})$ and
$\mathcal{E}|_F\cong\mathcal{O}_F(1)^{\oplus 4}$ for any fibre
$F\cong\mathbb{P}^{4}$ of $Y\to C$;

\item[$(9')$] $n=8$, there is a surjective morphism $q:Y\to\Gamma$
onto a smooth curve $\Gamma$ such that any general fibre $F$ of
$q$ is a smooth quadric hypersurface $\mathbb{Q}^{4}$ in
$\mathbb{P}^5$ with $\mathcal{E}|_F\cong\mathcal{O}_F(1)^{\oplus
4}$;

\item[$(10')$] $n=8$, there is a vector bundle $\mathcal{V}$ on a
smooth surface $S$ such that $Y\cong\mathbb{P}_S(\mathcal{V})$ and
$\mathcal{E}|_F\cong\mathcal{O}_F(1)^{\oplus 4}$ for any fibre
$F\cong\mathbb{P}^{3}$ of $Y\to S$;

\item[$(11')$] $n=8$, there exists a smooth projective $5$-fold
$W$ and a morphism $\pi:Y\to W$ expressing $Y$ as blown up at a
finite set $B$ of points and an ample vector bundle $\mathcal{E}'$
on $W$ such that $\mathcal{E} = \pi^*\mathcal{E}'\otimes
[-\pi^{-1}(B)]$ and $K_W+\det\mathcal{E}'$ is ample; or

\item[$(12')$] $n=7$ and $(Y,\mathcal{E})$ is as in
\em{Proposition} \ref{class-scroll k=m-3}.

\end{enumerate}
\end{Cor}

\noindent\textbf{Proof.} If $n\geq 11$, then $k=n-5\geq
6>\mathrm{dim } Y$, so $(X,L)$ is an adjunction scroll over $Y$ by
\cite[(2.1.1)]{BS2}. If $7\leq n\leq 10$, then we conclude by
\cite[(2.1.2), (2.1.3)]{BS2}, \cite[proposition 3.1]{T} and
Proposition \ref{class-scroll k=m-3}. $\square$

Finally, as to the case $n=6$ and $m=5$, under the extra
assumption that
$\mathrm{Pic}(Y)\cong\mathbb{Z}[\mathcal{O}_Y(1)]$, we can easily
deduce from Theorem \ref{thm} the following

\begin{Cor}\label{corollary}
Let $X$ be a smooth $6$-fold and $L$ an ample line bundle on $X$.
Assume that $(X,L)$ $\cong$ $(\mathbb{P}
(\mathcal{E}),\mathcal{O}_{\mathbb{P} (\mathcal{E})}(1))$ is a
$\mathbb{P}^1$-bundle, $\pi:X\to Y$, over a smooth $5$-fold $Y,
\mathcal{E}=\pi_*L$. If
$\mathrm{Pic}(Y)\cong\mathbb{Z}[\mathcal{O}_Y(1)]$, then $(X,L)$
is an adjunction-theoretic scroll over $Y$ under $\pi$ except
when:
\begin{enumerate}
\item[$(1'')$] $Y\cong\mathbb{P}^5$ and
$\mathcal{E}\cong\mathcal{O}_{\mathbb{P}^5}(a_1)\oplus\mathcal{O}_{\mathbb{P}^5}(a_2)$
with $a_1+a_2\leq 5$ and $a_i>0$ for $i=1,2$; \item[$(2'')$]
$Y\cong\mathbb{Q}^5$ and
$\mathcal{E}\cong\mathcal{O}_{\mathbb{Q}^5}(a_1)\oplus\mathcal{O}_{\mathbb{Q}^5}(a_2)$,
with $a_1+a_2\leq 4$ and $a_i>0$ for $i=1,2$, where $\mathbb{Q}^5$
is a smooth quadric hypersurface of $\mathbb{P}^6$; \item[$(3'')$]
$Y$ is a Del Pezzo $5$-fold, i.e. $-K_Y\simeq\mathcal{O}_Y(4)$,
and either $(\alpha)$ $\mathcal{E}\cong\mathcal{O}_Y(1)^{\oplus
2}$, or $(\beta)$
$\mathcal{E}|_{l}\cong\mathcal{O}_{\mathbb{P}^1}(1)\oplus\mathcal{O}_{\mathbb{P}^1}(2)$
for any line $l\cong\mathbb{P}^1$ of $(Y,\mathcal{O}_Y(1))$;
\item[$(4'')$] $Y$ is a Mukai $5$-fold, i.e.
$-K_Y\simeq\mathcal{O}_Y(3)$, and
$\mathcal{E}\cong\mathcal{O}_Y(1)^{\oplus 2}$; \item[$(5'')$] $Y$
is a Fano $5$-fold such that $-K_Y\simeq\det\mathcal{E}$ and $X$
is a Fano $6$-fold such that $-K_X\simeq 2L$.
\end{enumerate}
\end{Cor}

\noindent\textbf{Proof.} Since
$\mathrm{Pic}(Y)\cong\mathbb{Z}[\mathcal{O}_Y(1)]$, if $(X,L)$ is
not an adjunction-theoretic scroll over $Y$ under $\pi$, then from
the formula (\ref{a4}) we deduce that
$K_Y+\det\mathcal{E}=\mathcal{O}_Y(a)$ for some integer $a\leq 0$.
Therefore, we have either $a=0$, obtaining Case $(5'')$ of the
statement, or $a<0$. In the latter case, since
$K_Y+\det\mathcal{E}$ is not nef, we conclude by Theorem
\ref{thm}. $\square$

\begin{Rem}\label{remark}
If $(X,L)$ and $(Y, \mathcal{E})$ are as in {\em{Corollary
\ref{corollary}}}, then similar techniques and arguments used for
the proofs of {\em \cite[theorem 3.1]{BS2}} and {\em
\cite[proposition 3.4]{T}} work also in the case
$\mathrm{Pic}(Y)\not\cong\mathbb{Z}[\mathcal{O}_Y(1)]$. More
precisely, we can say in this situation when $(X,L)$ is not an
adjunction-theoretic scroll over $Y$ under the projection map
$\pi:X\to Y$, but the final classification result appears at
present not yet complete in some of its parts. This fact is due
especially to the lack of classification results of pairs
$(Y,\mathcal{E})$ as above with $K_Y+\det\mathcal{E}$ trivial, and
of properties on the third reduction of $(Y,\mathcal{E})$ in the
same spirit of {\em \cite[remark 4]{AN}}. We refer to {\em
\cite[$\S 3$]{AN}} and {\em \cite{Ohno}} for the case $m-2\geq 2$
and the second reduction of the pair $(Y,\mathcal{E})$.
\end{Rem}

\bigskip

\noindent\textbf{Acknowledgement.} I would thank Proff. M.C.
Beltrametti, A. Lanteri and A. Laface for their kind comments and
some useful remarks about the final form of this paper. Finally,
the author wants to thank for some 
funds supporting this research from the University of Milan (FIRST
2006 and 2007).


\bigskip

\bigskip


\bigskip
\bigskip

%

Departamento de Matem\'{a}tica

Universidad de Concepci\'{o}n

Casilla 160-C, Concepci\'{o}n (Chile)

\medskip

\textit{E-mail address:} atironi@udec.cl


\begin{thebibliography}{10}

\bibitem{APW} V. Ancona, T. Peternell, J.A.
Wi$\acute{\textrm{s}}$niewski. Fano bundles and splitting theorems
on projective spaces and quadrics. \textit{Pacific Journal of
Math.} \textbf{163} (1994), no. 1, 17--42.

\bibitem{AM} M. Andreatta, M. Mella. Contractions on a manifold polarized
by an ample vector bundle. \textit{Trans. Amer. Math. Soc.}
\textbf{349} (1997), no. 11, 4669--4683.

\bibitem{AN} M. Andreatta, C. Novelli. Manifolds polarized by
vector bundles. \textit{Ann. Mat. Pura Appl. (4)} \textbf{186}
(2007), no. 2, 281--288.

\bibitem{AO} M. Andreatta, G. Occhetta. Ample vector bundles with sections
vanishing on special varieties. \textit{Internat. J. Math.}
\textbf{10} (1999), no. 6, 677--696.

\bibitem{AW} M. Andreatta, J.A. Wi$\acute{\textrm{s}}$niewski.
On manifolds whose tangent bundle contains an ample subbundle.
\textit{Invent. Math.} \textbf{146} (2001), no. 1, 209--217.

\bibitem{BS} M.C. Beltrametti, A.J. Sommese. The Adjunction
Theory of Complex Projective Varieties. Expositions in
Mathematics, vol. 16, W. de Gruyter, Berlin--New York (1995).

\bibitem{BS2} M.C. Beltrametti, A.J. Sommese. Comparing the
classical and the adjunction theoretic definition of scrolls,
Geometry of complex projective varieties (Cetraro, 1990), 55--74.

\bibitem{BSW} M.C. Beltrametti, A.J. Sommese, J.A.
Wi$\acute{\textrm{s}}$niewski. Results on varieties with many
lines and their applications to adjunction theory (with an
appendix by M.C. Beltrametti and A.J. Sommese), in Complex
Algebraic Varieties, Bayreuth 1990, ed. by K. Hulek, T. Peternell,
M. Schneider, and F.-O. Schreyer, \textit{Lecture Notes in Math.}
\textbf{1507} (1992), 16--38.

\bibitem{CMS} K. Cho, Y. Miyaoka, N.I. Shepherd-Barron. Characterizations
of projective space and applications to complex symplectic
manifolds. Higher dimensional birational geometry (Kyoto, 1997),
Math. Soc. Japan, Tokyo, \textit{Adv. Stud. Pure Math.}
\textbf{35} (2002), 1--88.

\bibitem{Fujita} T. Fujita. On polarized manifolds whose adjoint bundles
are not semipositive. Algebraic Geometry, Sendai, 1985,
\textit{Adv. Stud. Pure Math.} \textbf{10} (1987), 167--178.

\bibitem{F1} T. Fujita. On adjoint bundles of ample vector bundles. In Proc.
Alg. Geom. Conf. Bayreuth, \textit{Lect. Notes Math.}
\textbf{1507} (1990), 105--112.

\bibitem{F2} T. Fujita. Classification Theories of Polarized Varieties,
\textit{London Math. Soc. Lecture Notes Series} \textbf{155},
Cambridge Univ. Press, Cambridge (1990)

\bibitem{I} P. Ionescu. Generalized adjunction and applications. \textit{Math.
Proc. Cambridge Philos. Soc.} \textbf{99} (1986), no. 3, 457--472.

\bibitem{KS} Y. Kachi, E. Sato. Segre's reflexivity and an
inductive characterization of hyperquadrics. \textit{Mem. Amer.
Math. Soc.} \textbf{160} (2002), n. 763.

\bibitem{LM} A. Lanteri, H. Maeda. Special varieties in adjunction theory
and ample vector bundles. \textit{Math. Proc. Cambridge Philos.
Soc.} \textbf{130} (2001), 61--75.

\bibitem{LM1} A. Lanteri, H. Maeda. Ample vector bundle characterizations
of projective bundles and quadric fibrations over curves.
\textit{Higher Dimensional Complex Varieties, Trento, 1994}, (ed.
M. Andreatta \and T. Peternell), (Walter de Gruyter, 1996), pp.
247--259.

\bibitem{Maeda} H. Maeda. Nefness of adjoint bundles for ample vector
bundles. \textit{Le Matematiche} Vol. L, Fasc. I (1995), 73--82.

\bibitem{Mella} M. Mella. Existence of good divisors on Mukai varieties.
\textit{J. Algebraic Geom.} \textbf{8} (1999), no. 2, 197--206.

\bibitem{Miyaoka} Y. Miyaoka. Numerical characterisations of hyperquadrics. Complex
analysis in several variables--Memorial Conference of Kiyoshi
Oka's Centennial Birthday (Tokyo) (Math. Soc. Japan, Ed.),
North-Holland, \textit{Adv. Stud. Pure Math.} \textbf{42} (2004),
209--235.

\bibitem{Miyaoka1} Y. Miyaoka. Numerical characterisations of
hyperquadrics. The Fano Conference, Univ. Torino, Turin (2004),
559--562.

\bibitem{Mori} S. Mori. Threefolds whose canonical bundles are not
numerically effective. \textit{Ann. of Math. (2)} \textbf{116}
(1982), no. 1, 133--176.

\bibitem{Mori1} S. Mori. Threefolds whose canonical bundles are not numerically
effective, in Algebraic Threefolds, Proceedings Varenna, 1981, ed.
by A. Conte, \textit{Lecture notes in Math.} \textbf{947} (1982),
125--189, Springer-Verlag, New York.

\bibitem{NO} C. Novelli, G. Occhetta. Ruled Fano fivefolds of index two.
\textit{Indiana Univ. Math. J.} \textbf{56} (2007), no. 1,
207--241.


\bibitem{Occhetta} G. Occhetta. A note on the classification of Fano manifolds of
middle index. \textit{Manuscripta Math.} \textbf{117} (2005), no.
1, 43--49.

\bibitem{Occhetta1} G. Occhetta. On some Fano manifolds of large pseudoindex.
\textit{Manuscripta Math.} \textbf{104} (2001), 111--121.

\bibitem{Ohno} M. Ohno. Classification of generalized polarized manifolds by
their nef values. \textit{Adv. Geom.} \textbf{6} (2006), no. 4,
543--599; Addendum to: ``Classification of generalized polarized
manifolds by their nef values'' [\textit{Adv. Geom.} \textbf{6}
(2006), no. 4, 543--599], \textit{Adv. Geom.} \textbf{7} (2007),
no. 2, 315--316.

\bibitem{OSS} C. Okonek, M. Schneider, H. Spindler. Vector bundles on
complex projective spaces. \textit{Progress in Mathematics}
\textbf{3} (1980), Birkhäuser, Boston, Mass.

\bibitem{Ot} G. Ottaviani. Spinor bundles on quadrics. \textit{Trans. Amer.
Math. Soc.} \textbf{307} (1988), no. 1, 301--316.

\bibitem{P} T. Peternell. A characterization of $\mathbb{P}^n$ by vector
bundles. \textit{Math. Z.} \textbf{205} (1990), no. 3, 487--490.

\bibitem{P1} T. Peternell. Ample vector bundles on Fano manifolds.
\textit{Internat. J. Math.} \textbf{2} (1991), no. 3, 311--322.

\bibitem{PSW} T. Peternell, M. Szurek, J.A. Wi$\acute{\textrm{s}}$niewski.
Fano manifolds and vector bundles. \textit{Math. Ann.}
\textbf{294} (1992), 151--165.

\bibitem{PSW1} T. Peternell, M. Szurek, J.A. Wi$\acute{\textrm{s}}$niewski.
Numerically effective vector bundles with small Chern classes.
Complex algebraic varieties (Bayreuth, 1990), \textit{Lecture
Notes in Math.} \textbf{1507}, 145--156, Springer, Berlin, 1992.

\bibitem{Pr} A. Predonzan. Intorno agli $S_k$ giacenti sulla varietà
intersezione completa di pi\`{u} forme. \textit{Atti Accad. Naz.
Lincei. Rend. Cl. Sci. Fis. Mat. Nat.} (8) \textbf{5} (1948),
238--242.

\bibitem{Sommese} A.J. Sommese. On the adjunction theoretic structure of projective
varieties. In Complex Analysis and Algebraic Geometry, Proceedings
G$\ddot{\textrm{o}}$ttingen 1985, \textit{Lecture Notes in Math.}
\textbf{1194} (1986), 175--213.

\bibitem{S} A.J. Sommese. On the adjunction theoretic structure of projective
varieties, Complex Analysis and Algebraic Geometry, Proceedings
G$\ddot{\textrm{o}}$ttingen 1985, Springer--Verlag,
\textit{Lecture Notes in Math.} \textbf{1194} (1986), 175--213.

\bibitem{T} A.L. Tironi. Scrolls over a four dimensional variety.
\textit{Adv. Geom.} \textbf{10} (2010), 145--159.

\bibitem{W} J.A. Wi$\acute{\textrm{s}}$niewski. Length of extremal rays and
generalized adjunction. \textit{Math. Z.} \textbf{200} (1989),
409--427.

\bibitem{W1} J.A. Wi$\acute{\textrm{s}}$niewski. On a conjecture of
Mukai. \textit{Manuscripta Math.} \textbf{68} (1990), 135--141.

\bibitem{W2} J.A. Wi$\acute{\textrm{s}}$niewski. On Fano manifolds of large
index. \textit{Manuscripta Math.} \textbf{70} (1991), 145--152.

\bibitem{YZ} Y.-G. Ye, Q. Zhang. On ample vector bundles whose adjunction
bundles are not numerically effective. \textit{Duke Math. J.}
\textbf{60} (1990), no. 3, 671--687.

\bibitem{Zhang} Q. Zhang. A theorem on the adjoint system for vector
bundles. \textit{Manuscripta Math.} \textbf{70} (1991), no. 2,
189--201.

\end{thebibliography}
\end{document}